\documentclass[aap,MSNbibl,nameyear,dvips]{arximspdf}

\doi{10.1214/13-AAP923} %kopijuoti is PTS
\volume{24}
\issue{2}
\pubyear{2014}
\firstpage{553}
\lastpage{598}

\makeatletter
\newcommand{\rrVert}{\Vert}
\newcommand{\rrvert}{\vert}
\newcommand{\llVert}{\Vert}
\newcommand{\llvert}{\vert}
\renewcommand{\citep}[1]{[\citet{#1}]}
\renewcommand{\citeyearpar}[1]{(\citeyear{#1})}
\newcommand{\eqref}[1]{(\ref{#1})}
\newtheorem{theorem}{Theorem}[section]
\newtheorem{lemma}[theorem]{Lemma}
\newtheorem{corollary}[theorem]{Corollary}
\newtheorem{proposition}[theorem]{Proposition}
\newproclaim{remark}[theorem]{Remark}
\newtheorem{hypothesis}[theorem]{Hypothesis}
\newtheorem{standing}[theorem]{Standing Hypothesis}
\newproclaim{definition}[theorem]{Definition}
\def\R#1{{\mathbb R}^{#1}} %bold superscripted R-variable
\def\R{\mathbb{R}}
\makeatother

\begin{document}
\begin{frontmatter}

\title{Quasi-stationary distributions for randomly perturbed dynamical
systems}
\runtitle{QSD for random perturbations}

\begin{aug}
\author[A]{\fnms{Mathieu} \snm{Faure}\ead[label=e1]{mathieu.faure@univ-amu.fr}}
\and
\author[B]{\fnms{Sebastian J.} \snm{Schreiber}\corref{}\ead[label=e2]{sschreiber@ucdavis.edu}\thanksref{t1}}
\runauthor{M. Faure and S. J. Schreiber}
\affiliation{Aix-Marseille University (Aix-Marseille School of Economics),\break  CNRS \& EHESS and University of California, Davis}
\address[A]{GREQAM, centre de la vieille charit\'e\\
2 rue de la vieille charit\'e\\
13236 Marseille Cedex 02\\
France\\
\printead{e1}} %adresu isvedimo komanda gale!
\address[B]{Department of Evolution and Ecology\\
%Graduate Group in Applied Mathematics\\
University of California\\
Davis, California 95616\\
USA\\
\printead{e2}}
\end{aug}
\thankstext{t1}{Supported in part by NSF Grants EF-09-28987 and DMS-10-22639 and a
start-up grant from the Office of the Dean of the College of Biological Sciences at University of California, Davis.}

% HISTORY:
\received{\smonth{2} \syear{2012}}
\revised{\smonth{10} \syear{2012}}

% ABSTRACT
%
\begin{abstract}
We analyze quasi-stationary distributions
$\{\mu^\varepsilon\}_{\varepsilon>0}$ of a family of Markov chains
$\{X^\varepsilon\}_{\varepsilon>0}$ that are random perturbations of a
bounded, continuous map $F\dvtx M\to M$, where $M$ is a closed subset
of $\mathbb{R}^k$. Consistent with many models in biology, these Markov
chains have a closed absorbing set $M_0\subset M$ such that
$F(M_0)=M_0$ and $F(M\setminus M_0)=M\setminus M_0$. Under some large
deviations assumptions on the random perturbations, we show that, if
there exists a positive attractor for $F$ (i.e., an attractor for $F$
in $M\setminus M_0$), then the weak* limit points of $\mu_\varepsilon$
are supported by the positive attractors of $F$. To illustrate the
broad applicability of these results, we apply them to nonlinear
branching process models of metapopulations, competing species,
host-parasitoid interactions and evolutionary games.
\end{abstract}

% KEYWORDS

\end{frontmatter}

%s1 #&#
%s1 ###
\section{Introduction}\label{s1}
A fundamental issue in biology is what are the minimal conditions to
ensure the long-term survivorship for all of the interacting
components, whether they be viral particles, bio-chemicals, plants or
animals. When these conditions are met the interacting populations are
said to persist or coexist. Since the pioneering work of
\citet{Lotka-25} and \citet{volterra-26} on competitive and
predator-prey interactions, \citet{thompson-24},
\citet{nicholson-bailey-35} on host-parasite interactions and
\citet{kermack-mckendrick-27} on disease outbreaks, nonlinear
difference and differential equations have been used to understand
conditions for persistence of interacting populations. For these
deterministic models, persistence is often equated with an attractor
bounded away from the extinction states in which case persistence holds
over an infinite time horizon~\citep{jtb-06}. In reality, extinction in
finite time is inevitable due to finite population sizes and mortality
events occurring with positive probability. However, for systems with a
large number of individuals or particles, these times to extinction may
be sufficiently long so that the system remains in a ``metastable
state,'' bounded away from extinction for a long time. To provide a
rigorous mathematical basis for this intuition, we develop a general
theory for randomly perturbed dynamical systems with absorbing states.
Under the appropriate assumptions about the random perturbations, we
show that the existence of a positive attractor (i.e., an attractor
which is bounded away from extinction states) for the unperturbed
system implies two things as the number of individuals or particles
gets large. First, when they exist, quasi-stationary distributions
concentrate on the positive attractors of the unperturbed system.
Second, the expected time to extinction for systems starting according
to this quasi-stationary distribution grows exponentially with the
system size. In particular, we generalize earlier related work for
one-dimensional randomly perturbed dynamical
systems~\citep{hognas-97,klebaner-lazar-zeitouni-98,ramanan-zeitouni-99}
to higher dimensional systems by extending a general theory of randomly
perturbed systems without absorbing states
[Kifer~(\citeyear{kif88,kif89,kif90})] to a general theory of randomly
perturbed systems with absorbing states.

For the unperturbed, deterministic dynamics, we consider a bounded
continuous map $F\dvtx M \rightarrow M$, where $M$ is a closed subset
of $\mathbb{R}^d$. \textit{A random perturbation of} $F$ is a family of
homogeneous Markov chains $\{X^{\varepsilon}\}_{\varepsilon> 0}$ on
$M$, whose transition kernels
\[
p^{\varepsilon}(x,\Gamma) = \mathbb{P} \bigl[X_{n+1}^{\varepsilon}
\in \Gamma\mid X_n^{\varepsilon} = x \bigr],\qquad x \in M, \Gamma \mbox{ Borel subset of } M
\]
enjoy the following standing hypothesis.

%
%st1.1 #&#
\begin{standing} \label{H}
For any $\delta>0$, $\lim_{\varepsilon\to0}
\beta_\delta(\varepsilon)=0$ where
\[
\beta_{\delta}(\varepsilon) = \sup_{x \in M}
p^{\varepsilon} \bigl(x, M\setminus N^{\delta}\bigl(F(x)\bigr) \bigr)
\]
and $N^{\delta}(A):=  \{x \in M \dvtx \inf_{y \in A} \llVert x-y\rrVert
< \delta \}$ denotes the $\delta$-neighborhood of $A$.
\end{standing}

All standing hypotheses are assumed to hold throughout the paper.
Standing Hypothesis~\ref{H} implies that $p^{\varepsilon}(x,\cdot)$
converges uniformly to the Dirac measure $\delta_{F(x)}$ at $F(x)$ for
the vague convergence of measures, that is, for any continuous function
$g\dvtx M\to\R$ with compact support,
\[
\lim_{\varepsilon\to0}\sup_{x \in M} \biggl\llvert \int
_M g(y) p^{\varepsilon
}(x,dy) - g\bigl(F(x)\bigr)\biggr\rrvert =0.
\]
Consequently, for small $\varepsilon>0$, the asymptotic behavior of the
Markov chain $\{X_n^{\varepsilon}\}_{n=1}^\infty$ should be related to
the dynamics of iterating the map $F$.

When the perturbed system admits an invariant measure (e.g., the Markov
chains are irreducible), the correspondence between the deterministic
dynamics and the perturbed counterpart was initially studied by
\citet{PAV33}, and more recently by
\citet{Sin72,Rue81,FreWen84} and Kifer
(\citeyear{kif88,kif89,kif90}). Recall that a Borel probability measure
$\mu_{\varepsilon}$ on $M$ is called \textit{a stationary distribution
for} $p^{\varepsilon}$ if
\[
\int_M p^{\varepsilon}(x,\Gamma) \mu_{\varepsilon}(dx)
= \mu _{\varepsilon} (\Gamma)\qquad\mbox{for any Borel set }\Gamma\subset M.
\]
These invariant measures describe the long-term statistical behavior of
the perturbed system. Let us assume, for a moment, that for all
$\varepsilon>0$, the Markov chain $X^{\varepsilon}$ admits (at least)
one invariant measure $\mu_{\varepsilon}$. We call $\mu$ a
\textit{limiting measure} of the family of measure
$\{\mu_{\varepsilon}\}_{\varepsilon>0}$ if $\mu$ is the weak* limit of
some sequence $\{\mu_{\varepsilon_n}\}_{n=1}^\infty$, where
$\varepsilon_n$ decreases to zero. Natural questions about these
limiting measures include: Are the limiting measures invariant for the
deterministic dynamics? If so, what can be said about their support?

Kifer~(\citeyear{kif88,kif89,kif90}) considered these questions under
the assumption that the transition kernels $p^{\varepsilon}$ satisfied
the following large deviation assumption: there exists a continuous,
nonnegative rate function $\rho$ such that
%
%e1 #&#
%e1 ###
\begin{equation}
\label{LD} \lim_{\varepsilon\rightarrow0} \varepsilon\log p^{\varepsilon}(x,U) =
- \inf_{y \in U} \rho(x,y)
\end{equation}
for any open set $U \subset M$ and uniformly in $x \in M$. Under
suitable assumptions, Kifer proved that limiting measures are invariant
for $F$ [i.e., $\mu(\Gamma) = \mu(F^{-1}(\Gamma))$ for all Borel set
$\Gamma$] and are supported by the attractors of the deterministic
dynamics. In particular, Kifer's approach allowed him to derive some of
Freidlin and Wentzell's results on the asymptotic behavior of invariant
measures for diffusion processes $X_t^{\varepsilon}$ generated by
operators of the form $L^{\varepsilon} = \varepsilon L + b$ where $L$
is a ``good'' second-order elliptic differential operator and $b$ a
vector field [\citet{FreWen84}, Chapter~6].

While Kifer's results are applicable to a wide class of stochastic
models for the physical sciences, they are not applicable to many
models in ecology, epidemiology, immunology and evolution. These
stochastic models often have absorbing states $M_0\subset M$
corresponding to the loss of one or more populations that satisfy the
following standing hypothesis.

%
%st1.2 #&#
\begin{standing} \label{hy.M}
The state space $M$ can be written $M =M_0 \cup M_1$, where:
\begin{itemize}
\item $M_0$ is a closed subset of $M$;

\item $M_0$ and $M_1$ are positively $F$-invariant, that is,
    $F(M_0) \subseteq M_0$ and\break  $F(M_1) \subseteq M_1$;\vadjust{\goodbreak}

\item[$\bullet$] the set $M_0$ is assumed to be absorbing for the
random perturbations
%
%e2 #&#
%e2 ###
\begin{equation}
\label{hy.abs} p^{\varepsilon}(x,M_1)= 0\qquad\mbox{for all }\varepsilon>0, x \in M_0.
\end{equation}
\end{itemize}
\end{standing}

For many of these biological models, the set $M_0$ of absorbing states
is reached in finite time almost surely for any $\varepsilon>0$.
Despite this eventual absorption, the process
$\{X_n^\varepsilon\}_{n=1}^\infty$ may spend an exceptionally long
period of time in the set $M_1$ of transient states provided that
$\varepsilon>0$ is sufficiently small. In applications, this
``metastable'' behavior may correspond to long-term persistence of an
endemic disease, coexistence of interacting species, or maintenance of
a genetic polymorphism. One approach to examining metastable behavior
are quasi-stationary distributions which are invariant distributions
when the perturbed process is conditioned on nonabsorption. More
specifically, we have the following.

%
%de1.3 #&#
\begin{definition}\label{df.qsd}
A probability measure $\mu_{\varepsilon}$ on $M_1$ is a
\textit{quasi-stationary distribution} (QSD) for $p^{\varepsilon}$
provided there exists $\lambda_{\varepsilon}\in(0,1) $ such that
\[
\int_{M_1} p^{\varepsilon}(x,\Gamma) \mu_{\varepsilon}(dx)
= \lambda _{\varepsilon} \mu_{\varepsilon}(\Gamma)\qquad\mbox{for all Borel }\Gamma \subset M_1.
\]
\end{definition}

Equivalently, dropping the index $\varepsilon$, a QSD $\mu$ satisfies
the identity
\[
\mu(\Gamma) = \mathbb{P}_{\mu} (X_n \in\Gamma\mid X_n \in M_1 )\qquad\forall n,
\]
where $\mathbb{P}_{\mu}$ denotes the law of the Markov chain
$\{X_n\}_{n=0}^\infty$, conditional to $X_0$ being distributed
according to $\mu$. Quasi-stationary distributions can sometimes be
defined through the so-called \textit{Yaglom limit},
\[
\mu(\Gamma) = \lim_{n \rightarrow+ \infty} \mathbb{P}_x
(X_n \in \Gamma\mid X_n \in M_1 ),
\]
when the limit exists and is independent of the initial state $x \in
M_1$. When \mbox{$\mathbb{P}(X_1 \in\cdot) = \mu(\cdot)$}, $\lambda$ is
the probability of not being absorbed in the next time step. The
existence of QSDs has been studied extensively
[\citet{DarSen65,SenVer66,Tweedie1974,Barbour1976,Nummelin1976,Arjas1980,Kijima1992,Feretal95,Chan1998,LasPea01,Gosselin2001,coolen2006quasi,Buckley2010limit}].

\citet{hognas-97,klebaner-lazar-zeitouni-98,ramanan-zeitouni-99}
studied weak* limit points $\mu$ of QSDs $\mu_\varepsilon$ as
$\varepsilon\to0$ for maps of the interval, that is, $M=[0,1]$ and
$M_0=\{0\}$. Under suitable assumptions, these authors have shown that
if $F$ admits an attractor in $(0,1)$, then the limiting measure $\mu$
is \mbox{$F$-}invariant and concentrated on the attractors of $F$ in
$(0,1)$. Moreover, $\lambda_\varepsilon\geq 1-e^{-c/\varepsilon}$ for
an appropriate constant $c>0$. This final assertion implies that if the
perturbed processes is initiated in the quasi-stationary state, then
the expected time to absorption increases exponentially with exponent
$1/\varepsilon$ as $\varepsilon$ decreases to zero.

Here, we extend these types of results to higher dimensional systems
where $M$~is a subset of $\R^d$. The two main results of the paper are
stated in Section~\ref{s2}. First, we state a general result that
ensures that the QSDs concentrate on attractors of $F$ restricted to
$M_1$. This result requires conditions on the topological dynamics and
the rate at which $\beta_\delta(\varepsilon)$ in Standing
Hypothesis~\ref{H} goes to zero. Second, for many applications, the
randomly perturbed Markov chains satisfy large deviation assumptions.
We present a result that guarantees the conditions of the general
theorem are satisfied. Proofs of these two results are presented in
Sections~\ref{s3} and~\ref{s4}, respectively. In Section~\ref{s3}, we
also show how the main result of
\citet{klebaner-lazar-zeitouni-98} can be derived from our general
theorem. In Section~\ref{s5}, we apply our results to stochastic models
of metapopulation dynamics, competing species, host-parasitoid
interactions and evolutionary games. In Section~\ref{s6}, we conclude
by verifying the large deviation assumptions for the examples in
Section~\ref{s5}.

%%%%%%%%%%%%%%%%%%%%%%%%%%%%%%%%%%%%%%%%%%%%%%%%%%%%%%%%%%%%%%%%
%%%%%%%%%%%%%%%%%%%%%%%%%%%%%%%%%%%%%%%%%%%%%%%%%%%%%%%%%%%%%%%%
%s2 #&#
%s2 ###
\section{Statement of the main results}\label{s2}

Let $\{X^{\varepsilon}\}_{\varepsilon>0}$ be a family of Markov chains
on a closed set $M \subset\R^d$, which satisfies
%whose transition kernels $p^{\varepsilon}$ satisfy Hypotheses \ref{H}
%and
Standing Hypothesis \ref{hy.M}. Since $M$ is assumed to be a closed
subset of $\R^d$, every topological concept must be understood in terms
of the topology induced in $M$. In particular, in the following, a
compact set $K$ will always be a closed (in $M$) bounded set $K \subset
M$.
%As a consequence, if $K$ is a compact set (in $M$) such that $K

We assume that, for each $\varepsilon>0$, there exists at least one QSD
$\mu_{\varepsilon}$: there exists $1 > \lambda_{\varepsilon} > 0$ such
that $\lambda_\varepsilon\mu_{\varepsilon}=Q^\varepsilon
\mu_\varepsilon$ where $Q^{\varepsilon}$ is the operator defined on
the set of finite Borel measures on $M_1$ by
\[
Q^{\varepsilon}(\mu) (\Gamma) = \int_{M_1}p^{\varepsilon}(x,
\Gamma) \mu (dx)\qquad\mbox{for every Borel } \Gamma\subset M_1.
\]
%
%where $q^{\varepsilon}$ is the restriction of the transition kernel
%$p^{\varepsilon}$ to $M_1$: $\forall x \in M_1, A_1 \in

Our main results concern characterizing the support of weak* limit
points $\mu$ of the $\mu_\varepsilon$ as $\varepsilon\downarrow0$.
Under suitable assumptions, we show that these weak* limit points are
supported by attractors of the map $F$ that lie in $M_1$; see
Section~\ref{s2.1} for a definition of an attractor. In
Section~\ref{s2.1}, we describe sufficient conditions for this result
with suitable assumptions about the topological dynamics of $F$ and
$\beta_\delta(\varepsilon)$ in Standing Hypothesis~\ref{H} goes to
zero. In Section~\ref{s2.2}, we describe large deviation assumptions which
satisfy the assumptions presented Section~\ref{s2.1} and which are
easier to verify for applications presented in Section~\ref{s5}.

%%%%%%%%%%%%%%%%%%%%%%%%%%%%%%%%%%%%%%%%%%%%%%%%%%%%%%%%%%%%%%%%%%%%%%%%%%%%%%%%%%%%%%%%%%%%%%%%%%%%%
%s2.1 #&#
%s2.1 ###
\subsection{Absorption-preserving chain recurrence and convergence to
attractors}\label{s2.1}

We begin by recalling a few definitions from dynamical system theory.
Let $F^n$ be the $n$-iterate of $F$. A set $B\subset M$ is
\textit{invariant} for $F$ if $F(B)=B$. A compact set $A$ is an
\textit{attractor for $F$} provided there exists an open neighborhood
$U$ of $A$ such that $\bigcap_{n\geq 1}F^n(U)=A$ and, for any open set
$V \supset A$, there exists $n(V)$ such that $F^n(U) \subset V$ for all
$n \geq n(V)$.

The key notions needed for our main result is \textit{absorption
preserving pseudoorbits and chain recurrence} introduced in
\citet{jacSch06}. These definitions generalize Conley's
\citeyearpar{Con78} notion of pseudoorbits and chain recurrence. Given
$\delta>0$, a family of points $\xi=(\xi_0,\ldots,\xi_{n}) \in M^{n+1}$
is called an \textit{absorption preserving $\delta$-pseudoorbit joining
$x$ to $y$} (ap $\delta$-pseudoorbit for short) provided that:
\begin{longlist}[(a)]
\item[(a)] $x= \xi_0, y = \xi_{n}$,

\item[(b)] $\xi_i \in M_0 \Rightarrow\xi_{i+1} \in M_0$ and

\item[(c)] $d(\xi_{i+1},F(\xi_i)) < \delta$, $i=0,\ldots,n-1$.
\end{longlist}
One can think of ap $\delta$-pseudoorbits as approximations of actual
orbits of the dynamics of $F$ with an error no greater than $\delta$
and that preserve the absorbing set~$M_0$. For readers unfamiliar with
these concepts, consider $F$ to be the identity map on the interval
$[0,1]$. Then any two points on the interval are connected by
\mbox{$\delta$-}pseudoorbits, for any $\delta>0$. However, as every
point is a fixed point, none of the points are connected by interating
the map $F$.

Given $x,y \in M$, we say that $x$ \textit{ap-chains to} $y$ and write
$x <_{\mathrm{ap}} y$ if for any $\delta>0$, there exists an ap
$\delta$-pseudoorbit joining $x$ to $y$. Notice that no point in $M_0$
ap-chains to any point in $M_1$. If $x <_{\mathrm{ap}} y$ and $y
<_{\mathrm{ap}} x$, we shall write $x \sim_{\mathrm{ap}} y$. If $x
\sim_{\mathrm{ap}} x$, then $x$ is an \textit{ap-chain recurrent}
point. The set $\mathcal{R}_{\mathrm{ap}}$ of ap-chain recurrent points
is $F$-invariant.The relation $\sim_{\mathrm{ap}}$, restricted to this
set defines an equivalence relation. The equivalent classes,
$[x]_{\mathrm{ap}}$ with $x\in\mathcal{R}_{\mathrm{ap}}$, are called
\textit{ap-basic classes}. In Section~\ref{s3.1}, we prove various
properties of these equivalence classes, for example,
$\omega(x)\subset\mathcal{R}_{\mathrm{ap}}$ whenever $\omega(x)\subset
M_0$ or $\omega(x) \subset M_1$.

Let $[x]_{\mathrm{ap}}$ and $[y]_{\mathrm{ap}}$ be two distinct
ap-basic classes. We write $[x]_{\mathrm{ap}} <_{\mathrm{ap}}
[y]_{\mathrm{ap}}$ if $x <_{\mathrm{ap}} y$. A maximal basic class
$[x]_{\mathrm{ap}}$ (i.e., $[x]_{\mathrm{ap}} <_{\mathrm{ap}}
[y]_{\mathrm{ap}}$ implies that $[x]_{\mathrm{ap}} =
[y]_{\mathrm{ap}}$) is called an \textit{ap-quasiattractor}. In
general, an ap-quasiattractor need not be an attractor for $F$. A
simple example is an increasing function $F\dvtx[0,1]\to[0,1]$ with
$F(x)=x$ for $x=1-1/n$ for all natural numbers $n$ and $F(x)\neq x$
otherwise. If $M_0=\{0\}$, then $x=1$ is a quasi-attractor but not an
attractor for $F$.

We need three hypotheses in order to state the first main result. The
first hypothesis requires that there is a finite number of ap-basic
classes including at least one ap-quasiattractor in $M_1$. This
assumption is satisfied for many important classes of mappings,
including gradient-like systems and Axiom A systems. When this
hypothesis is satisfied, we prove in Section~\ref{s3.2} that all the
ap-quasiattractors are in fact attractors.

%
%hy2.1 #&#
\begin{hypothesis} \label{hy.finite} There exists only a finite number
of ap-basic classes in $M_1 \dvtx \{K_i\}_{i=1,\ldots,v}$. Moreover, we
assume that they are closed sets and $\{K_i\}_{i=1,\ldots,\ell}$ with
$\ell\geq 1$ are the ap-quasiattractors and
$\{K_i\}_{i=\ell+1,\ldots,v}$ are the nonap-quasi\-attractors.
\end{hypothesis}

Our second hypothesis ensures the time spent near nonap-quasiattractors
is not too long relative the $\beta_\delta(\varepsilon)$ described in
Standing Hypothesis~\ref{H}. For any Borel set~$V$, we define the first
passage time $\tau^\varepsilon_V= \min\{ n\dvtx X_n^\varepsilon \notin
V\}$.

%
%hy2.2 #&#
\begin{hypothesis} \label{hy.na} Given any $\delta>0$, there exist
neighborhoods $V_{i}\subset N^\delta(K_i)$ of $K_i$ for $\ell+1\leq  i
\leq  v$ and quantity $\delta_1\in(0,1)$ such that
\[
\sup_{x \in V_j} \mathbb{P}_x \bigl[
\tau^\varepsilon_{V_j} > h(\varepsilon) \bigr] \leq\zeta(
\varepsilon) \quad\mbox{and}\quad\lim_{\varepsilon\to0} \zeta(\varepsilon)=0
\]
for a function $h$ satisfying
\[
\lim_{\varepsilon\to0 }h(\varepsilon) \beta_{\delta_1}(\varepsilon)
=0.
\]
\end{hypothesis}

Our final hypothesis provides a lower bound on the probability of
absorption on the event $X_n^\varepsilon$ is sufficiently close to
$M_0$.

%
%hy2.3 #&#
\begin{hypothesis} \label{hy.boundary} There exists a neighborhood
$V_0$ of $M_0$ such that
\[
\lim_{\varepsilon\to0} \frac{\beta_{\delta_0}(\varepsilon)}{\inf_{x
\in V_0} p^{\varepsilon}(x,M_0)} = 0.
\]
\end{hypothesis}

We prove in Section~\ref{s3} that, if $M_0$ is a global attractor, then
$\mu$ is supported by~$M_0$; see Theorem \ref{th.M0}. The main result
of this section is the following theorem. A proof is given in
Section~\ref{s3}.

%
%th2.4 #&#
\begin{theorem} \label{th.suppmu} Assume that Hypotheses \ref{hy.finite} and \ref{hy.na} hold. Then any weak* limit point $\mu$ of
$\{\mu^\varepsilon\}_{\varepsilon>0}$ satifies $\mu(V_j) = 0$ for
$j=\ell+1,\ldots,v$. In addition, if Hypothesis~\ref{hy.boundary}
holds, then $\mu$ is supported by the union of the attractors
$\bigcup_{i=1}^\ell K_i$. Moreover, there exists a $\delta>0$ such that
$\lambda _\varepsilon\geq 1-\beta_\delta(\varepsilon)$ for all
$\varepsilon>0$ sufficiently small.
\end{theorem}

%s2.2 #&#
%s2.2 ###
\subsection{Large deviation hypotheses}\label{s2.2}
For applications, it is often easier to verify certain large deviation
hypotheses rather than Hypotheses~\ref{hy.na} and \ref{hy.boundary}. To
this end we consider the following large deviation hypothesis.

%
%hy2.5 #&#
\begin{hypothesis} \label{hy.rho}
There exists a function $\rho\dvtx M
\times M \rightarrow[0,+\infty]$ such that:
\begin{longlist}[(iii)]
\item[(i)] $\rho$ is continuous on $M_1 \times M$,

\item[(ii)] $\rho(x,y) = 0$ if and only if $y = F(x)$,\vadjust{\goodbreak}

\item[(iii)] for any $\beta>0$,
%
%e3 #&#
%e3 ###
\begin{equation}
\label{beta} \inf \bigl\{ \rho(x,y) \dvtx x \in M, y \in M, d\bigl(F(x),y\bigr) >
\beta \bigr\} >0,
\end{equation}
where $d(x,y)=\max_i |x_i-y_i|$,

\item[(iv)] for any open set $U$, we have the large deviations lower
    bound
%
%e4 #&#
%e4 ###
\begin{equation}
\label{lower} \liminf_{\varepsilon\rightarrow0} \varepsilon\log
p^{\varepsilon}(x,U) \geq- \inf_{y \in U} \rho(x,y)
\end{equation}
that holds uniformly for $x$ in compact subsets of $M_1$
whenever $U$ is an open ball in $M$. Additionally, for any closed
set $C$, we have the uniform upper bound
%
%e5 #&#
%e5 ###
\begin{equation}
\label{upper} \limsup_{\varepsilon\rightarrow0} \sup_{x \in M}
\varepsilon\log p^{\varepsilon}(x,C) \leq- \inf_{y \in C}
\rho(x,y).
\end{equation}
\end{longlist}
\end{hypothesis}

Equations \eqref{beta} and \eqref{upper} imply that Standing
Hypothesis \ref{H} holds.
Additionally, since $M_0$ is absorbing, \eqref{lower} implies that
$\rho(x,y) = + \infty$ for all $x \in M_0$, $y \in M_1$. The upper
bound \eqref{upper} can be weakened as a uniform bound on compact
subsets of $M_1$. In that case, Hypothesis \ref{H} is no longer implied
by Hypothesis \ref{hy.rho}.

We also make the following assumption that ensures absorption is
reasonably likely when the process is near the absorbing states.

%
%hy2.6 #&#
\begin{hypothesis} \label{hy.rho2}
For any $c>0$, there exists an open neighborhood $V_0$ of~$M_0$ such
that
%
%e6 #&#
%e6 ###
\begin{equation}
\label{V_0} \lim_{\varepsilon\rightarrow0} \inf
_{x \in V_0} \varepsilon\log p^{\varepsilon}(x,M_0)
\geq-c.
\end{equation}
\end{hypothesis}

To state our main result under these large deviation assumptions, we
need to introduce an alternative notion of chain recurrence. Given $n
\in\mathbb{N}^*=\{0,1,2,3,\ldots\}$, define the function $A_n$ on
$M^{n+1}=\underbrace{M\times\cdots\times M}_{n+1\ \mathrm{times}}$ by
\[
\xi= (\xi_0,\ldots,\xi_{n}) \mapsto A_n(\xi)
= \sum_{i=0}^{n-1} \rho(\xi_i,
\xi_{i+1}).
\]
$A_n$ measures the ``cost'' of
$X^\varepsilon_n$ following the partial trajectory $\xi$ where the cost
is measured in terms of how much ``noise'' is required to move along
this partial trajectory. For any $x,y$ in $M$, we define
\[
B_{\rho}(x,y) = \inf \bigl\{A_n(\xi) \mid n \geq1, \xi\in
M^{n+1}, \xi _0 = x, \xi_n = y \bigr\}.
\]
The function $B_\rho(x,y)$ represents the minimal cost in going from
$x$ to $y$. $B_\rho$~induces a partial order on $M$ by writing $x
<_{\rho} y$ (i.e., ``$x$ $\rho$-chains to $y$'') if \mbox{$B_{\rho}(x,y)
= 0$}. Roughly, $x$ $\rho$-chains to $y$ if the,re exist paths joining
$x$ to $y$ with arbitrarily low costs. If $x <_{\rho} y$ and $y
<_{\rho} x$, we write $x \sim_{\rho} y$.

We define the \textit{set of $\rho$-chain recurrent points
$\mathcal{R}_{\rho}$} to be the set of points $x \in M$ such that $x
\sim_{\rho} x$. The \textit{$\rho$-basic classes} are the equivalence
classes for $\sim_{\rho}$ restricted to the $\rho$-chain recurrent set.
Since a point in $M_0$ never $\rho$-chains to a point in $M_1$, the
$\rho$-basic classes are included either in $M_0$ or in $M_1$. In
general, the $\rho$-basic classes and the ap-basic classes introduced
in Section~\ref{s2.1} need not be equivalent. For example, consider
$F\dvtx [0,1]\to[0,1]$ given by the identity map $F(x)=x$ for all $x$
and $M_0=\varnothing$. Let $\rho(x,y)=|x-y|$. Then every point $\{x\}$
is a $\rho$-basic class, but the only ap-basic class is $[0,1]$.
However, unlike this example, if there is a finite number of
$\rho$-basic classes, then we prove in Section~\ref{s4} (see Theorem
\ref{th.equalitybasic}) that the ap-basic classes and $\rho$-basic
classes agree.

Given a $\rho$-chain recurrent point $x$, let $[x]_\rho$ denote its
$\rho$-basic class. We say that $[x]_{\rho} <_{\rho} [y]_{\rho}$ if $x
<_{\rho} y$ and call \textit{$\rho$-quasiattractors} the maximal
$\rho$-equivalence classes. When a $\rho$-quasiattractor $A$ is
isolated (i.e., there is a neighborhood of the quasi-attractor
containing no other $\rho$-chain recurrent point), we prove in
Section~\ref{s4} that $A$ is an \textit{attractor} for $F$; see
Proposition~\ref{pr.rhoQA=A}.

In Section~\ref{s4}, we use Theorem~\ref{th.suppmu} to prove the
following result. Applications of Theorem~\ref{th.ld} are given in
Section~\ref{s5}.

%
%th2.7 #&#
\begin{theorem} \label{th.ld} Assume that Hypotheses \ref{hy.rho} and
\ref{hy.rho2} hold and that there exists a finite number of $\rho
$-basic classes in $M_1$, which are closed. If:
\begin{itemize}
\item there is at least one $\rho$-quasiattractor $A$ among the
    $\rho$-basic classes in $M_1$, and

\item $\mu_\varepsilon(U)>0$ for any neighborhood $U$ of $A$ and
$\varepsilon>0$,
\end{itemize}
then any weak*-limit point of $\{\mu_{\varepsilon}\}_{\varepsilon>0}$
is $F$-invariant and is supported by the union of $\rho
$-quasiattractors in $M_1$. Moreover, there exists $c>0$ such that
$\lambda_\varepsilon\geq 1- e^{-c/\varepsilon}$ for all
$\varepsilon>0$.
\end{theorem}

%
%re2.8 #&#
\begin{remark}\label{remark.ld}
Assume that there is a finite number of closed nonquasiattractors
$[x_1]_{\rho},\ldots,[x_N]_{\rho}$ in $M_1$ and $A=(\mathcal{R}_{\rho}
\cap M_1)\setminus U_{i=1}^N [x_i]_{\rho}$ is an attractor for $F$.
Then the main result still holds: if $\mu _\varepsilon(U)>0$ for any
neighborhood $U$ of $A$ and $\varepsilon>0$, then any weak*-limit point
of $\{\mu_{\varepsilon}\}_{\varepsilon>0}$ is $F$-invariant and is
supported by~$A$. Moreover, there exists $C>0$ such that
$\lambda_\varepsilon\geq 1- e^{-C/\varepsilon}$ for all
$\varepsilon>0$.
\end{remark}

%%%%%%%%%%%%%%%%%%%%%%%%%%%%%%%%%%%%%%%%%%%%%%%%%%%%%%%%%%%%%%%%%
%%%%%%%%%%%%%%%%%%%%%%%%%%%%%%%%%%%%%%%%%%%%%%%%%%%%%%%%%%%%%%%%%
%s3 #&#
%s3 ###
\section{\texorpdfstring{Proof of Theorem~\protect\ref{th.suppmu}}
{Proof of Theorem 2.4}}\label{s3}

In this section, we prove Theorem~\ref{th.suppmu}. We begin by proving
several key results about ap-chain recurrence in Sections~\ref{s3.1}
and~\ref{s3.2}. In Section~\ref{s3.3}, we prove some key properties of
limiting quasi-stationary distributions. A proof of
Theorem~\ref{th.suppmu} is given in Section~\ref{s3.4}. In
Section~\ref{s3.5}, we show how our proof of Theorem~\ref{th.suppmu}
provides an alternative proof of the main result of
\citet{klebaner-lazar-zeitouni-98}. In addition to the Standing
Hypotheses, the results in Section~\ref{s3.2} require
Hypothesis~\ref{hy.finite}, and the proofs in Sections~\ref{s3.4}
and~\ref{s3.5} require Hypotheses~\ref{hy.finite}, \ref{hy.na} and
\ref{hy.boundary}.

%%%%%%%%%%%%%%%%%%%%%%%%%%%%%%%%%%%%%%%%%%%%%%%%%%%%%%%
%s3.1 #&#
%s3.1 ###
\subsection{Absorption preserving chain recurrence}\label{s3.1}
We recall a few definitions and facts from dynamical systems. The
\textit{$\omega$-limit set} of $B \subset M$ is given by $\omega(B) =
\bigcap_{n \geq1} \overline{\bigcup_{p \geq n} F^p(B)}$. It is the
maximal invariant set in the closure of $\bigcup_{n \geq0}F^n(B)$. An
equivalent definition of an attractor presented in Section~\ref{s2.1}
is that a compact set $A$ is an attractor for $F$ provided it admits an
open neighborhood $U$ such that $\omega(U) = A$; the open set $\{x \in
M \dvtx \omega(x) \subset A\}$ is then called the \textit{basin of
attraction} of $A$. By a classical result [see \citet{Con78}], a
compact set $A$ is an attractor for $F$ if and only if there exists an
open set $V$ which contains $A$ and such that
%
%e7 #&#
%e7 ###
\begin{equation}
\label{attractor} F(\overline{V}) \subset V,\qquad \bigcap_{n \in\mathbb{N}}
F^n(V) = A.
\end{equation}

Our assumption that $\|F\|=\sup_{x\in M} \|F(x)\|<\infty$ implies that
the set $\mathcal{R}_{\mathrm{ap}}$ of ap-chain recurrent points is
included in $N^{\|F\|}(0)$. The relation $\sim_{\mathrm{ap}}$,
restricted to this set defines an equivalence relation. Unlike the
ap-basic classes lying in $M_0$, the ap-basic classes lying in $M_1$
may not be closed. However, we have the following:

%
%le3.1 #&#
\begin{lemma} \label{lm.apbasic=closed} Let $x$ be an ap-chain
recurrent point in $M_1$. Then $\overline{[x]_{\mathrm{ap}}} \subset
M_0 \cup [x]_{\mathrm{ap}}$. In particular,
\[
\overline{[x]_{\mathrm{ap}}} \subset M_1 \Rightarrow[x]_{\mathrm{ap}}
\mbox{ closed}.
\]
\end{lemma}

\begin{pf}
Let $y \in\overline{[x]_{\mathrm{ap}}}$. There exists a sequence
$\{y_k\}$ in $[x]_{\mathrm{ap}}$ which converges to $y$. Any ap
$\delta$-pseudoorbit from $x$ to $y_k$ is an ap $2\delta$-pseudoorbit
from $x$ to $y$, provided $k$ is chosen large enough. Hence $x
<_{\mathrm{ap}} y$. On the other hand, assume that $y \notin M_0$ and
consider an ap $\delta$-pseudoorbit $(\xi_0,\ldots,\xi_{n})$ chaining
$y_k$ to $x$. We have
\[
d\bigl(F(y),\xi_1\bigr) \leq\delta+ d\bigl(F(y),F(y_k)
\bigr) \leq2\delta
\]
by continuity of $F$ provided $k$ is large enough. Consequently,
$(y,\xi_1,\ldots,\xi_n)$ is an ap $2\delta$-pseudoorbit chaining $y$ to
$x$ and $y\in[x]_{\mathrm{ap}}$.
\end{pf}

The following lemma shows that ap-basic classes are invariant.

%
%le3.2 #&#
\begin{lemma} Any ap-basic class $[x]_{\mathrm{ap}}$ is positively
$F$-invariant: $F([x]_{\mathrm{ap}}) \subseteq[x]_{\mathrm{ap}}$. If
$\overline {[x]_{\mathrm{ap}}} \subset M_1$ (which implies that
$[x]_{\mathrm{ap}}$ is closed), it is $F$-invariant,
$F([x]_{\mathrm{ap}})= [x]_{\mathrm{ap}}$.
\end{lemma}

\begin{pf} If $[x]_{\mathrm{ap}}$ is a singleton, then $F(x)=x$, and there is nothing
to prove. Assume that there exists $y \neq x$ such that $y \in
[x]_{\mathrm{ap}}$. For any $\delta>0$, continuity and boundedness of
$F(M)$ imply that there exists a $\delta/2>\delta'>0$ such that
\[
d\bigl(z,F(x)\bigr) < \delta' \Rightarrow d\bigl(F(z),F^2(x)
\bigr)<\delta/2\qquad\mbox{for all }z\in M.
\]
Pick an ap $\delta'$-pseudoorbit $(x=\xi_0,\xi_{1},\ldots,\xi_n=y)$
joining $x$ to $y$. Since\break  $d(\xi_1, F(x))\leq \delta'$
\[
d\bigl(F^2(x),\xi_2\bigr) \leq d\bigl(F^2(x),F(
\xi_1)\bigr)+ d\bigl(F(\xi_1),\xi_2\bigr)
\leq\delta /2 + \delta' < \delta
\]
and $(F(x),\xi_2,\ldots,\xi_n)$ is an ap $\delta$-pseudoorbit joining
$F(x)$ to $y$. Hence,\break  $F([x]_{\mathrm{ap}})
\subseteq[x]_{\mathrm{ap}}$.

Next, let us assume that $[x]_{\mathrm{ap}}$ is closed in $M_1$. For
every $y\in [x]_{\mathrm{ap}}$, we need to prove that $y= F(y')$ for
some $y' \in [x]_{\mathrm{ap}}$. For any $\delta>0$, choose an ap
$\delta$-pseudoorbit $(\xi^{\delta}_i)_{i=0,\ldots,n(\delta)}$ joining
$y$ to itself. Now choose a compact set $K \subset M_1$ containing an
open neighborhood of $[x]_{\mathrm{ap}}$. We prove in the next section
(see Remark \ref{rk.ap}) that the families $\xi^{\delta}$ can be chosen
in such a way that they are contained in $K$. In particular the family
$\xi^{\delta}_{n(\delta) -1}$ admits an accumulation point $y'$ as
$\delta\to0$. By continuity of $F$, $F(y') = y$ and, therefore, $y'
\sim_{\mathrm{ap}} y\sim_{\mathrm{ap}} x$.
\end{pf}

For classical chain recurrence, $\omega(x)$ is contained in the chain
recurrent set. While ap-chain recurrence shares this property whenever
$\omega(x)\subset M_0$ or $\omega(x) \subset M_1$, in general it only
satisfies a weaker property.

%
%le3.3 #&#
\begin{lemma}\label{lem.omega} For $x\in M$, $\omega(x)\cap\mathcal
{R}_{\mathrm{ap}}\neq\varnothing$.
\end{lemma}

\begin{pf} If $x\in M_0$ or $\omega(x)\subset M_1$, then the classical result
for chain recurrence implies
$\omega(x)\subset\mathcal{R}_{\mathrm{ap}}$. Suppose $x\in M_1$ and
$y\in\omega(x)\cap M_0$. Then $\omega(y)
\subset\mathcal{R}_{\mathrm{ap}}$. Since $\omega(y)\subseteq\omega(x)$,
the result follows.
\end{pf}

%Let $[x]_{\mathrm{ap}}$ and $[y]_{\mathrm{ap}}$ be two distinct ap-basic classes. We
%write $[x]_{\mathrm{ap}} <_{\mathrm{ap}} [y]_{\mathrm{ap}}$ if $x <_{\mathrm{ap}} y$. A maximal basic
%class $[x]_{\mathrm{ap}}$ (i.e. $[x]_{\mathrm{ap}} <_{\mathrm{ap}} [y]_{\mathrm{ap}}$ implies that
%$[x]_{\mathrm{ap}} = [y]_{\mathrm{ap}}$) is called an \textit{ap-quasiattractor}.

%
%le3.4 #&#
\begin{lemma} If $[x]_{\mathrm{ap}}$ is maximal, then $x <_{\mathrm{ap}} z$ if and only
if $z\in[x]_{\mathrm{ap}}$. In particular, any ap-quasiattractor
$[x]_{\mathrm{ap}}$ is compact.
\end{lemma}

\begin{pf} Let $z$ be such that $x <_{\mathrm{ap}} z$. To prove that $z \in[x]_{\mathrm{ap}}$,
we need to show that $z<_{\mathrm{ap}}x$. By Lemma \ref{lem.omega},
$\omega(z) \cap\mathcal{R}_{\mathrm{ap}} \neq\varnothing$. Hence there
exists $z' \in \omega(z) \cap\mathcal{R}_{\mathrm{ap}}$. In particular,
$x <_{\mathrm{ap}} z <_{\mathrm{ap}} z'$. As $z' \in\mathcal{R
}_{\mathrm{ap}}$, maximality of $[x]_{\mathrm{ap}}$ implies that $z'
\in[x]_{\mathrm{ap}}$. Thus $z <_{\mathrm{ap}} x$. In particular, if
$y\in \overline{[x]_{\mathrm{ap}}}$, then the proof of
Lemma~\ref{lm.boundary} implies that $y>_{\mathrm{ap}} x$. Hence,
$y\in[x]_{\mathrm{ap}}$ and $[x]_{\mathrm{ap}}$ is closed.
\end{pf}

%[x]_{\mathrm{ap}}$. Consequently $z \in[x]_{\mathrm{ap}}$. Otherwise, let $z' \in M_0
%$[x]_{\mathrm{ap}} <_{\mathrm{ap}} \omega(z')$, we once again have $z \in[x]_{\mathrm{ap}}$.

The next result is an easy consequence of Proposition~4.2 in
\citet{kif88}. A~closed ap basic set is said to be isolated in
$M_1$ if it admits an open neighborhood which is disjoint from any
other ap basic class:

%
%th3.5 #&#
\begin{theorem} \label{th.qa=a} Let $[x]_{\mathrm{ap}}$ be an isolated
ap-quasiattractor in $M_1$. Then $[x]_{\mathrm{ap}}$ is an attractor.
\end{theorem}

%%%%%%%%%%%%%%%%%%%%%%%%%%%%%%%%%%%%%%%%%%%%%%%%%%%%%%%%%%%%%%%
%s3.2 #&#
%s3.2 ###
\subsection{Finiteness of the ap-basic classes}\label{s3.2}
Throughout this subsection, we require Hypothesis~\ref{hy.finite}.
Namely, there exists a finite number of ap-basic classes
$\{K_i\}_{i=1}^v$ where the $K_i$ are closed sets, $\{K_i\}_{i=1}^\ell$
are ap-quasiattractors and $\{K_i\}_{i=\ell+1}^v$ are
nonap-quasiattractors.

The following result is proved in \citeauthor{kif88}
[(\citeyear{kif88}), pages 217--218] for classical chain recurrence. We
give a proof adapted to our settings for the convenience of the reader.

%
%le3.6 #&#
\begin{lemma} \label{lm.217}
\textup{(a)} For any $\theta>0$ sufficiently small, there exists a
    quantity $0<\delta(\theta)<\theta$ such that, if there is an ap
    $\delta(\theta)$-pseudoorbit $(\xi_0,\ldots,\xi_{n})$
    satisfying
%
%e8 #&#
\begin{eqnarray}\label{delta3} %
d(\xi_0,K_i) &<& \delta(\theta),\qquad d(\xi_j,K_i) >\theta\quad\mbox{and}
\nonumber\\[-10pt]\\[-10pt]
d(\xi_{n}, K_{i'}) &<& \delta(\theta)\qquad\mbox{for some }i,i' \in\{1,\ldots,v\}, j \in\{1,\ldots,n\},\nonumber
\end{eqnarray}
then $i\neq i'$ and $K_{i'} >_{\mathrm{ap}} K_i$.

\textup{(b)} For any $\delta'>0$, there exist $\delta\in(0,
    \delta')$ and $n_0\geq 1$ such that any ap $\delta$-pseudo\-orbit
    of length greater than $n_0$ must pass through
    $N^{\delta'}(\mathcal{R}_{\mathrm{ap}})$.
\end{lemma}

\begin{pf}
Assume that, for any $\delta>0$, there exists an ap
$\delta$-pseudoorbit $(\xi^{\delta}_0,\ldots,\xi^{\delta}_{n(\delta)})$
such that
%
%e9 #&#
%e8 ###
\begin{equation}
\label{delta2} d\bigl(\xi^{\delta}_0,K_i\bigr)
\leq\delta \quad\mbox{and}\quad d\bigl(\xi^{\delta}_{n(\delta
)},K_{i'}
\bigr) \leq\delta.
\end{equation}
Then there exists $\delta_k \downarrow0$, $y \in K_i$ and $y' \in
K_{i'}$ such that $\lim_{k\to\infty}\xi^{\delta_k}_0 =y$
and~$\lim_{k\to\infty}\xi_{n(\delta_k)}^{\delta_k} = y'$. Hence
$d(F(y),\xi^{\delta_k}_1) < \delta_k + d(F(y),F(\xi_0^{\delta_k}))$
and\break  $d(F(\xi^{\delta_k}_{n(\delta_k)-1}),y') \leq \delta_k +
d(\xi_{n(\delta_k)}^{\delta_k},y')$. Therefore, for any $\delta>0$,
$(y,\xi^{\delta_k}_1,\ldots,\break  \xi_{n(\delta_k)-1}^{\delta_k},y')$
is an ap $\delta$-pseudoorbit provided that $k$ is large enough. This
proves that $K_{i'} >_{\mathrm{ap}} K_i$. As a consequence, if $K_{i'}
>_{\mathrm{ap}} K_i $ does not hold, this means that there exists some
quantity $\overline{\delta}>0$ such that, for any $0<\delta
<\overline{\delta}$, we cant have (\ref{delta2}). Now pick $\theta>0$
smaller than $\overline{\delta}$.

Now assume that $i=i'$. Choose $\theta$ small enough such that
$\overline{N^{\theta}(K_i)} \subset M_1$. Assume that there exist a
decreasing sequence $\delta_k \downarrow0$ and ap
$\delta_k$-pseudoorbits $(\xi^k_{0},\ldots,\xi^k_{n_k})$ such that
(\ref{delta3}) holds with $\xi= \xi^k$, $\delta= \delta_k$ and $j =
j_k$. Without loss of generality, we may assume that
$\lim_{k\to\infty}\xi^k_0 = x \in K_i$, $\lim_{k\to\infty}\xi^k_{j_k}
=y \in K \setminus N^{\theta}(K_i)$ and $\lim_{k\to\infty}\xi^k_{n_k} =
z \in K_{i}$, where $K\subset M_1$ is a compact set such that
$F(N^\theta(K_i))\subsetneq K$. We then have $d(\xi^k_1,F(x))
\leq\delta_k + d(F(\xi_0^k),F(x))$, $d(y, F(\xi^k_{j_k-1})) \leq
d(y,\xi^k_{j_k}) + \delta_k$, $d(F(y),\xi^k_{j_k+1}) \leq
d(F(y),F(\xi^k_{j_k})) + \delta_k$ and\break  $d(z, F(\xi^k_{n_k-1})) \leq
d(y,\xi^k_{j_k}) + \delta_k$. By continuity of $F$, this implies that,
for any $\delta>0$, the sequence
$(x,\xi^k_1,\ldots,\xi^k_{j_k-1},y,\xi^k_{j_k+1},\ldots,\xi^k_{n_k-1},z)$
is an\vspace*{1pt} ap $\delta$-pseudoorbit, provided $k$ is large enough.
Consequently, $x<_{\mathrm{ap}}y<_{\mathrm{ap}}z$ contradicting the
fact that $K_i$ is an ap-basic class.

We now prove point (b). For any $x \in M$ and $\gamma>0$,
Lemma~\ref{lem.omega} implies that the quantity
\[
n^{\gamma}(x):= \inf \bigl\{n \in\mathbb{N} \dvtx F^n(x) \in
N^{\gamma
}(\mathcal{R}_{\mathrm{ap}}) \bigr\}
\]
is finite. By continuity of $F$, $n^\gamma$ is upper-semicontinuous.
Compactness of $\overline{F(M)}$ and upper semicontinuity imply that
\[
n^{\gamma}:= \max_{x\in M} n^\gamma(x) \leq\max
_{y \in
\overline{F(M)}} n^{\gamma}(y) + 1 < \infty.
\]

Now assume that there exists $\delta'>0$ such that the statement of
(b) is not true. In particular, for each $k$ there exists an ap
$\delta_k=\delta'/k$-pseudoorbit of length $n^{\delta'/2}$, $\xi^k =
(\xi_0^k,\ldots,\xi^k_{n^{\delta'/2}})$, which does not enter
$N^{\delta'}(\mathcal{R}_{\mathrm{ap}}) $. Passing\vspace*{-3pt} to a
subsequence if necessary, we may assume that $\lim_{k\to\infty} \xi^k_j
= \xi_j\in M$ for any $j=1,\ldots,n^{\delta'/2}$. The sequence $\xi$ is
a partial solution of the discrete dynamical system induced by $F$,
that is, $F(\xi_i)=\xi_{i+1}$ for $i=0,\ldots,n^{\delta'/2}-1$. The
definition of $n^{\delta'/2}$ implies that there exists $j_0$ such that
$d(\xi_{j_0}, \mathcal{R}_{\mathrm{ap}}) \leq\delta'/2$. Hence,
$\xi^k_{j_0} \in N^{\delta'}(\mathcal{R}_{\mathrm{ap}})$ for $k$ large
enough which contradicts the choice of $\xi^k$.
\end{pf}

%
%re3.7 #&#
\begin{remark}\label{rk.ap} Notice that, even without the
finiteness assumption, the following statement still holds: given an
ap-basic class $[x]_{\mathrm{ap}}$ in $M_1$ and $\theta>0$, there
exists a quantity $\delta>0$ such that any ap $\delta$-pseudoorbit
joining $[x]_{\mathrm{ap}}$ to itself remains into
$N^{\theta}([x]_{\mathrm{ap}})$.
\end{remark}

%
%co3.8 #&#
\begin{corollary} \label{summary} Given $\delta'>0$, there exist
isolating open neighborhoods $\{V_i\}_{i=1,\ldots,v}$ of the ap-basic
classes $\{K_i\}_{i=1,\ldots,v}$, and positive constants $\delta_1$ and
$n_0$ such that:
\begin{longlist}[(a)]
\item[(a)] $\overline{N^{\delta_1}(K_i)}\subset V_i$ for $1\leq
    i\leq  v$;

\item[(b)] any ap $\delta_1$-pseudoorbit starting in $V_i$
    remains in $V_i$ for $i=1,\ldots,\ell$;

\item[(c)] if there exists an ap $\delta_1$-pseudoorbit
    $(\xi_0,\ldots,\xi_{n})$ such that $\xi_0 \in
    N^{\delta_1}(K_i)$, $\xi_{n} \in N^{\delta_1}(K_{i'})$ and
    $\xi_k \notin V_i$ for some $2\leq  k\leq  n-1$, then $i \neq
    i'$ and $K_{i'} > K_i$.

\item[(d)] any\vspace*{1pt} ap $\delta_1$-pseudoorbit of length greater than
    $n_0$ must pass through
    $N^{\delta'}(\mathcal{R}_{\mathrm{ap}})$.
\end{longlist}
\end{corollary}

\begin{pf} Choose $\theta\in(0,\delta')$ sufficiently small so that
Lemma~\ref{lm.217}(a) holds, and let $\delta(\theta)>0$ be as given by
Lemma~\ref{lm.217}(a). Choose neighborhoods $V_i$ of $K_i$ such that
$\overline{N^\theta(K_i)}\subset V_i$ for $i=1,\ldots,v$ and
$F(\overline{V_i})\subset V_i$ for $i=1,\ldots,k$. The latter choice is
possible as Lemma~\ref{lm.apbasic=closed} implies that the ap-basic
sets $K_i$ are compact for $i=1,\ldots,v$, and Theorem~\ref{th.qa=a}
implies that $K_i$ is an attractor for $i=1,\ldots,\ell$. Choose
$\delta_1\in(0,\delta(\theta))$ such that $\delta_1$ is less than the
$\delta$ given by Lemma~\ref{lm.217}(b) and such that any ap
$\delta_1$-pseudoorbit starting in $V_i$ for $i=1,\ldots,\ell$ remains
in $V_i$. This latter choice is possible as $F(\overline{V_i})\subset
V_i$ for $i=1,\ldots,\ell$.
\end{pf}

%%%%%%%%%%%%%%%%%%%%%%%%%%%%%%%%%%%%%%%%%%%%%%%%%%%%%%%%%%%
%s3.3 #&#
%s3.3 ###
\subsection{Limit behavior of quasi-stationary distributions}\label{s3.3}

Throughout this section, we assume that there exists a decreasing
sequence $\varepsilon_n \downarrow0$ such that, for every $n \in
\mathbb{N}$, $\mu_n$ is a quasi stationary probability measure for
$p^{\varepsilon_n}$ with associated eigenvalue $\lambda_n$.
Additionally, we assume that $\mu_n$ converges weakly to a Borel
probability measure~$\mu$. We note that the results in this subsection
do not require Hypotheses \ref{hy.finite}~or~\ref{hy.boundary}. Recall from Standing Hypothesis~\ref{H}
that $\beta_{\delta}(\varepsilon) = \sup_{x \in M} p^{\varepsilon} (x,
M\setminus N^{\delta}(F(x)) )$.

%
%le3.9 #&#
\begin{lemma} \label{lm.lambda}
We have the following:
\begin{longlist}[(a)]
\item[(a)] $\liminf_{n\to\infty} \lambda_n \geq\mu(M_1)$. In
    particular, if $\mu$ is supported by $M_1$, then
    $\lim_{n\to\infty}\lambda_n =1$. Alternatively, if
    $\lim_{n\to\infty}\lambda_n =0$, then $\mu$ is supported by $M_0$.
\item[(b)] If there exists an attractor $A \subset M_1$ such that
    $\mu_n(U) >0$ for every $n$ and every open neighborhood $U$ of $A$,
    then there exists $\delta>0$ such that
\[
\lambda_n \geq1 - \beta_{\delta}(\varepsilon_n)
\]
for all $n$.

\item[(b$'$)] If, in addition to the assumption of \textup{(b)}, there exists
    some neighborhood $V_0$ of $M_0$ such that
\[
\lim_{n\to\infty}\frac{\beta_{\delta}(\varepsilon_n)}{\inf_{x \in V_0}
p^{\varepsilon_n} (x,M_0)}= 0,
\]
then $\mu(V_0) = 0$.
\end{longlist}
\end{lemma}

\begin{pf}
(a) Let $(\delta_k)_k$ be a positive sequence, decreasing to zero, and
define
\[
V_k:= \bigl\{x \in M_1\dvtx d(x,M_0) >
\delta_k \bigr\},\qquad U_k:= F^{-1}(V_k).
\]
Notice that $(U_k)_k$ is an increasing sequence of open sets such that
$U_k \subset M_1$ (by $F$-invariance of $M_0$) and $\bigcup_{k}(U_k) =
M_1$ (by closedness of $M_0$). We have
\begin{eqnarray*}
\lambda_n &\geq& \int_{U_k}
\mu_n(dx) p^{\varepsilon_n}(x,M_1)
\\
&\geq& \mu_n(U_k) \inf_{x \in U_k}
p^{\varepsilon_n}(x,M_1)
\\
&=& \mu_n(U_k) \Bigl(1- \sup_{x \in U_k}
p^{\varepsilon_n}(x,M_0)\Bigr).
\end{eqnarray*}
Since $F(U_k) \subset V_k$, we have $N^{\delta_k}(F(U_k)) \subset M_1$.
Thus
\[
\lambda_n \geq\mu_n(U_k) \Bigl(1 - \sup_{x \in U_k} p^{\varepsilon_n} \bigl(x, N^{\delta_k}\bigl(F(x)
\bigr)^c\bigr)\Bigr).
\]
By weak* convergence, the definition of $\lambda_n$ and Standing
Hypothesis~\ref{H},
\[
\liminf_n \lambda_n \geq\liminf
_n \mu_n(U_k) \geq
\mu(U_k)
\]
for all $k$. Point (a) follows since $\lim_k \mu(U_k) =
\mu(M_1)$.\vadjust{\goodbreak}

(b) Choose an open neighborhood $U$ of $A$ such that
$F(\overline{U}) \subset U$ and $\delta>0$ such that
$N^{\delta}(F(\overline{U})) \subset U$. We have
\[
\lambda_n \mu_n(U) \geq\mu_n(U) \Bigl(1 -
\sup_{x \in\overline{U}} p^{\varepsilon_n}\bigl(x,U^c\bigr)\Bigr).
\]
Since $p^{\varepsilon_n}(x,U^c) \leq
p^{\varepsilon_n}(x,(N^{\delta}(F(\overline{U})))^c)$, and
$\mu_n(U)>0$, we get that $\lambda_n \geq1 -
\beta_{\delta}(\varepsilon_n)$.

(b$'$). By assumption, we have
\begin{eqnarray*}
1 - \beta_{\delta}(\varepsilon_n) &\leq& \lambda_n
\\
& =& \int_M \bigl(1 - p^{\varepsilon_n}(x,M_0)
\bigr) \mu_n(dx)
\\
&\leq& \mu_n(M \setminus V_0) + \mu_n(V_0)
\Bigl(1 - \inf_{x \in V_0} p^{\varepsilon_n}(x,M_0)
\Bigr),
\end{eqnarray*}
which gives
\[
\mu_n(V_0) \leq\frac{\beta_{\delta}(\varepsilon_n)}{\inf_{x \in V_0}
p^{\varepsilon_n}(x,M_0)}.
\]
Since $V_0$ is open and $\lim_{n\to\infty} \mu_n=\mu$ in the weak*
topology, the result follows.
\end{pf}

%
%re3.10 #&#
\begin{remark} Notice that we actually have a better result, as the
quantity $\beta_{\delta}(\varepsilon_n)$ could be replaced by the
smallest quantity
\[
\sup_{x \in\overline{U}}p^{\varepsilon_n}\bigl(x,\bigl(N^{\delta}
\bigl(F(\overline{U})\bigr)\bigr)^c\bigr).
\]
\end{remark}

%
%pr3.11 #&#
\begin{proposition}\label{prop.inv} If $\lim_{n\to\infty}\lambda_n =1$,
then the probability measure $\mu$ is \mbox{$F$-}invariant. In
particular, $\mu $ is supported by the closure of
$\mathcal{R}_{\mathrm{ap}}$.
\end{proposition}

\begin{pf} It suffices to verify that
%
%e10 #&#
%e9 ###
\begin{equation}
\label{eq.inv} \int_{M} g(x) \mu(dx) = \int
_{M} g \bigl(F (x)\bigr) \mu(dx)
\end{equation}
for any bounded continuous function $g\dvtx M \to\R$. Uniform
continuity of $g$ on $N^{\|F\|+ \delta}(0)$ and Hypothesis \ref{H}
imply
\[
\lim_{n\to\infty} \sup_x \biggl\llvert \int
_{M} \bigl(g(y) - g\bigl(F(x)\bigr)\bigr)
p^{\varepsilon
_n}(x,dy) \biggr\rrvert =0.
\]
Therefore,
\begin{eqnarray*}
&& \biggl\llvert \int_{M} \bigl(g(x) - g \bigl( F (x)\bigr)
\bigr)\mu_n(dx)\biggr\rrvert
\\
&&\qquad = \biggl\llvert \int
_{M} \biggl(\lambda_n \int_{M}
g(y) p^{\varepsilon_n}(x,dy) - g\bigl(F(x)\bigr) \biggr) \mu_n(dx)
\biggr\rrvert
\\
&&\qquad \leq 2(1-\lambda_n) \|g\| + \biggl\llvert \int
_{M} \biggl(\int_{M} \bigl(g(y) - g
\bigl(F(x)\bigr)\bigr) p^{\varepsilon_n}(x,dy) \biggr) \mu_n(dx) \biggr
\rrvert
\\
&&\qquad \leq 2(1-\lambda_n) \|g\| + \sup_x \biggl
\llvert \int_{M} \bigl(g(y) - g\bigl(F(x)\bigr)\bigr)
p^{\varepsilon_n}(x,dy)\biggr\rrvert.
\end{eqnarray*}
Sending $n$ to infinity implies \eqref{eq.inv} for any continuous
bounded $g$. Hence, $\mu$ is $F$-invariant. $F$-invariance of $\mu$
implies that the support of $\mu$ is contained in the Birkhoff center
of $F$, that is, the closure of recurrent points of $F$, $\{x \in M
\dvtx x \in\omega(x)\}$, which is in turn included in the closure of
$\mathcal{R}_{\mathrm{ap}}$.
\end{pf}

The following theorem provides a sufficient condition for the support
of the limiting measure $\mu$ to lie on the absorbing set $M_0$.

%
%th3.12 #&#
\begin{theorem}\label{th.M0} Assume that $M_0$ is a global attractor.
Then $\mu$ is supported by $M_0$.
\end{theorem}

\begin{pf} If $\liminf_{n\to\infty}\lambda_n =0$, Lemma~\ref{lm.lambda}
implies that $\mu(M_0)=1$. Assume that
$\liminf_{n\to\infty}\lambda_n>0$, and let $c=\inf_n \lambda_n>0$.
Given $\alpha>0$, pick an open neighborhood $U$ of $M_0$, $\delta_1
>0$ and $n_0 \in\mathbb{N}$ such that $U \subset N^{\alpha}(M_0)$,
$F(\overline{U}) \subset U$, any ap $\delta_1$-pseudoorbit starting in
$U$ remains in $U$ and any ap $\delta_1$-pseudoorbit of length at
least $n_0$ eventually enters $U$; see Corollary \ref{summary}.

Let $\mathcal{E}_{n,k}$ be the event $
\{(X_j^{\varepsilon_n})_{j=0,\ldots,k}\mbox{ is an ap }
\delta_1\mbox{-pseudoorbit}\}$. Since a $\delta_1$-pseudoorbit of
length at least $n_0$ ends in $U$, we have
\begin{eqnarray*}
&& \mathbb{P}_x \bigl[X_k^{\varepsilon_n} \in U^c \bigr]
\\
&&\qquad \leq \mathbb {P}_x \bigl[ \mathcal{E}_{n,k}^c \bigr] + \mathbb{P}
\bigl[ \mathcal {E}_{n,k}\mbox{ and }X_k^{\varepsilon_n} \in
U^c \bigr]
\\
&&\qquad \leq \sum_{j=0}^{k-1}
\mathbb{P}_x \bigl(d \bigl(X_{j+1}^{\varepsilon
_n}, F\bigl(X_j^{\varepsilon_n}\bigr) \bigr) > \delta_1 \bigr)+ 0
\\
&&\qquad \leq k \beta_{\delta_1}(\varepsilon_n)
\end{eqnarray*}
for $k \geq n_0$ and $x \in M$. The last inequality follows from the
definition of $\beta_\delta(\varepsilon)$, and the second inequality
from the fact that the event $\mathcal{E}^c_{n,k}$ is included in the
union of the $k$ events $\{d  (X_{j+1}^{\varepsilon_n},
F(X_j^{\varepsilon_n}) ) > \delta_1\}$. By the definition of
$\mu_n$,
\[
\mu_{n}\bigl(U^c\bigr) \leq\frac{1}{\lambda_n^{n_0}} \int
_M \mu_n(dx) \mathbb {P}_x
\bigl[X_{n_0}^{\varepsilon_n} \in U^c \bigr] \leq
\frac{n_0\beta
_{\delta_1}(\varepsilon_n)}{c^{n_0}}
\]
and the last quantity goes to zero as $n$ goes to infinity. Since
$\alpha>0$ was arbitrary, $\mu(M_1)=0$.
\end{pf}

%
%re3.13 #&#
\begin{remark} The proof is not needed in the particular case where
$\lambda_n$ goes to one since $\mu$ is then $F$-invariant and the
Birkhoff center is contained in $M_0$.
\end{remark}

%%%%%%%%%%%%%%%%%%%%%%%%%%%%%%%%%%%%%%%%%%%%%%%%%%%%%%%%%%%%%%%%%%%%%%%%%%%%%%%%%%%%%%%%%%
%s3.4 #&#
%s3.4 ###
\subsection{\texorpdfstring{Proof of Theorem~\protect\ref{th.suppmu}}
{Proof of Theorem 2.4}}\label{s3.4}
We assume Hypotheses~\ref{hy.finite}, \ref{hy.na} and \ref{hy.boundary}
hold. Recall that, under the finiteness assumption, the
ap-quasiattractors $\{K_i\}_{i=1,\ldots,\ell}$ are actually attractors;
see Theorem \ref{th.qa=a}. Also, there exists $\delta_0>0$ such that
$\lambda_n \geq1 - \beta_{\delta_0}(\varepsilon_n)$; see Lemma
\ref{lm.lambda}(b). Let $\{V_i\}_{i=1,\ldots, v}$ and $\delta_1 \leq
\delta_0$ be chosen as in Corollary \ref{summary}. Given a Borel set
$V$ define $\tau^n_V = \inf \{j \geq0 \dvtx X^{\varepsilon_n}_j \notin
V  \}$.

Call $b = v -\ell$ the number of nonap-quasiattractors in $M_1$ and
$K=\bigcup_{i=1}^v K_i$. Choose sequences $\{m_n\}_{n\geq 1}$ and
$\{m'_n\}_{n\geq 1}$ such that
\[
\lim_{n\to\infty} \beta_{\delta_1}(\varepsilon_n)
m_n = 0,\qquad \lim_{n\to \infty} \frac{m'_n}{m_n} =0 \quad\mbox{and}\quad \lim_{n\to\infty} \frac
{h(\varepsilon_n)}{m'_n} =0.
\]

The presence of an attractor inside $M_1$ such that $\mu_n(U) >0$ for
any $n$, and any open neighborhood $U$ implies that $\lim_{n\to\infty}
\lambda_n=1$, by Lemma~\ref{lm.lambda}(b). Proposition~\ref{prop.inv}
implies that $\mu$ is $F$-invariant and supported by the closure of
$\mathcal{R}_{\mathrm{ap}}$.

Let us prove the first statement of Theorem~\ref{th.suppmu}. Let $j
\in\{\ell+1,\ldots,v\}$ be fixed. By definition of $\lambda_n$,
\[
\mu_n(V_j) = \frac{1}{\lambda_n^{r}} \int
_{x \in M} \mu_n(dx) \mathbb {P}_x
\bigl[X_r^{\varepsilon_n} \in V_j \bigr]\qquad\forall r \in
\mathbb{N}^*.
\]

For $i=\{1,\ldots,b\}$, call $t_n^i$ the integer $\lfloor
m_n/i\rfloor$. Let $\mathcal{E}_n$ and $\mathcal{E}'_n$ be the events,
respectively, defined by
\[
\mathcal{E}_n = \bigl\{\bigl(X_i^{\varepsilon_n}
\bigr)_{i=1,\ldots,m_n} \mbox{ is a } \mbox{$\delta_1$-pseudoorbit}
\bigr\}
\]
and
\[
\mathcal{E}'_n = \bigl\{ \forall i \in\{\ell+1,\ldots,v
\},\ \forall q \geq m'_n, X_p^{\varepsilon_n}
\in N^{\delta_1}(K_i) \Rightarrow X_{p+q}^{\varepsilon_n}
\notin N^{\delta_1}(K_i) \bigr\}.
\]
The set $\mathcal{E}'_n$ is the event ``after its first entry in any
$N^{\delta_1}(K_i)$, the Markov chain will have escaped from this set
after $m'_n$ steps and will never come back.''

On the event $\mathcal{E}_n \cap\mathcal{E}'_n$, the process
$(X^{\varepsilon_n}_1,\ldots, X^{\varepsilon}_{m_n})$ is an ap
$\delta_1$-pseudoorbit and therefore gets trapped in
$\bigcup_{i=1}^\ell V_i$ if it enters in this set. Corollary
\ref{summary} implies that it cannot spend more than $b$ blocks of
length at most $m'_n$ in $\bigcup_{i=\ell+1}^v N^{\delta_1}(K_i)$. In
particular, for $n$ large enough, $X^{\varepsilon_n}_{m_n}$ is in $V_j$
only if $X^{\varepsilon}_{t_n^i} \in(N^{\delta_1}(K))^c$ for some $i
\in \{1,\ldots,b\}$. Therefore,
\[
\mathbb{P}_x \bigl[\bigl\{X_{m_n}^{\varepsilon_n} \in
V_j\bigr\} \cap\mathcal {E}_n \cap\mathcal{E}'_n
\bigr] \leq\sum_{i=1}^b
\mathbb{P}_x \bigl[X^{\varepsilon_n}_{t_n^i} \notin
N^{\delta_1}(K) \bigr].
\]

On the other hand on the event $\mathcal{E}_n$, starting from
$N^{\delta_1}(K_i)$, the chain cannot enter back into $N^{\delta
_1}(K_i)$ once it exited $V_i$ (by Corollary~\ref{summary}(c)).
Hypothesis~\ref{hy.na} implies
\[
\mathbb{P}_x\bigl[\bigl(\mathcal{E}'_n
\bigr)^c \cap\mathcal{E}_n\bigr] \leq\sum
_{i=\ell
+1}^v \sup_{y \in V_i}
\mathbb{P}_y \bigl[\tau^n_{V_i} \geq
m'_n \bigr] \leq b \zeta(\varepsilon_n)
\]
as $m'_n > h(\varepsilon_n)$ for $n$ sufficiently large. Consequently,
\begin{eqnarray*}
&& \mathbb{P}_x \bigl[X_{m_n}^{\varepsilon_n} \in V_j \bigr]
\\
&&\qquad \leq \mathbb{P}_x\bigl[(
\mathcal{E}_n)^c\bigr] + \mathbb{P}_x \bigl[
\bigl(\mathcal{E}'_n\bigr)^c \cap
\mathcal{E}_n \bigr] + \sum_{i=1}^b
\mathbb{P}_x \bigl[X^{\varepsilon_n}_{t_n^i} \notin
N^{\delta_1}(K) \bigr]
\\
&&\qquad \leq m_n \beta_{\delta_1}(\varepsilon_n) + b
\zeta(\varepsilon_n) + \sum_{i=1}^b
\mathbb{P}_x \bigl[X^{\varepsilon_n}_{t_n^i} \notin
N^{\delta_1}(K) \bigr]
\end{eqnarray*}
for $n$ sufficiently large. Therefore we have, using the invariance
properties of $\mu_n$,
\begin{eqnarray*}
&& \int\mu_n(dx) \mathbb{P}_x \bigl[X_{m_n}^{\varepsilon_n}
\in V_j \bigr]
\\
&&\qquad \leq m_n\beta_{\delta_1}(\varepsilon_n) + b \zeta(\varepsilon_n) + \sum
_{i=1}^b \int\mu_n(dx)\mathbb{P}_x \bigl[X^{\varepsilon
_n}_{t_n^i} \notin N^{\delta_1}(K) \bigr]
\\
&&\qquad \leq m_n \beta_{\delta_1}(\varepsilon_n) + b
\zeta(\varepsilon_n)+ \sum_{i=1}^b \lambda_n^{t_n^i}
\mu_n\bigl(\bigl(N^{\delta_1}(K)\bigr)^c\bigr)
\\
&&\qquad \leq m_n \beta_{\delta_1}(\varepsilon_n) + b
\zeta(\varepsilon_n) + b \mu_n\bigl(\bigl(N^{\delta_1}(K)
\bigr)^c\bigr),
\end{eqnarray*}
which converges to $0$ as $n\to\infty$. By our choice of the sequence
$m_n$,
\[
\liminf_{n\to\infty} \lambda_n^{m_n} \geq
\liminf_{n\to\infty} \bigl(1-\beta _{\delta_1}(
\varepsilon_n)\bigr)^{m_n} = 1.
\]
Hence, $\lim_{n\to\infty} \lambda_n^{m_n}=1$ and
\[
\lim_{n\to\infty} \mu_n(V_j) = \lim
_{n\to\infty} \frac{1}{\lambda
_n^{m_n}}\int\mu_n(dx)
\mathbb{P}_x\bigl[X_{m_n}^{\varepsilon_n} \in
V_j\bigr] =0.
\]
The proof of the first statement is complete since $\mu(V_j) \leq
\liminf_n \mu_n(V_j)$.

Now, under Hypothesis \ref{hy.boundary}, Lemma~\ref{lm.lambda} implies
that the support of $\mu$ is contained in $M \setminus V_0$. In
particular, $\mu(\mathcal{R}_{\mathrm{ap}} \cap M_1) = 1$. Hence,
$\mu(\bigcup_{i=1}^v K_j)=1$ and the second statement of
Theorem~\ref{th.suppmu} follows.

%%%%%%%%%%%%%%%%%%%%%%%%%%%%%%%%%%%%%%%%%%%%%%%%%%%%%%%%%%%%%%%%%%%%%%%%%%%%%%%%%%%%%%%%%%%%%%
%s3.5 #&#
%s3.5 ###
\subsection{\texorpdfstring{A derivation of Theorem 3 of Klebaner, Lazar and Zeitouni
(\protect\citeyear{klebaner-lazar-zeitouni-98})}{A derivation of
Theorem 3 of Klebaner, Lazar and Zeitouni (1998)}}\label{s3.5}

We assume Hypotheses~\ref{hy.finite}, \ref{hy.na} and \ref{hy.boundary}
hold. We can obtain Theorem 3 of
\citet{klebaner-lazar-zeitouni-98} as a consequence of our
proof\vadjust{\goodbreak}
of Theorem~\ref{th.suppmu}. To see why, we describe how their
assumptions (A1)--(A6) imply our main assumptions. For the sake of
brevity, we do not state their assumptions here. Rather we refer the
interested reader to their article. Under their assumptions (A1)--(A6),
\citet{klebaner-lazar-zeitouni-98} prove nonconvergence to the
finite set of unstable equilibria for one-dimensional maps. Their
assumptions (A1) and (A2) guarantee the following, calling $M_0 =
\{0,1\}$ to fit our settings,\setcounter{footnote}{1}\footnote{In the quoted paper, $M_1 =
[0,1]$ and $M_0$ is $[0,1]^c$, but it does not change the problem as
they consider continuous state space; see
\citet{ramanan-zeitouni-99}.} so we have:
\begin{longlist}[(a)]
\item[(a)] $f(M_0) \subset M_0$ and $f(M_1) \subset M_1$;
\item[(b)] the absorbing state $\{0\}$, the unstable equilibria
    $(x_i^*)_{i=0,\ldots,k}$, and the stable equilibria
    $(s_i)_{i=1,\ldots,l}$ form a Morse decomposition for the
    dynamical system induced by $f$; hence there is a finite number
    of ap-basic classes (see Proposition~\ref{prop.dyn}) and the
    ap-quasiattractors in $M_1$ are the $(s_i)_{i=1,\ldots,l}$. Our
    Hypothesis \ref{hy.finite} is verified.
\end{longlist}

By their assumption (A4), we derive the uniform (in $x$) upper bound
\[
\exists\lambda_0 >0, C>0\qquad \mbox{such that } P_{\varepsilon}\bigl(x, \bigl(N^{\delta}
\bigl(f(x)\bigr)\bigr)^c\bigr) \leq C e^{-\lambda_0 \delta/\varepsilon}
\]
for any $\delta>0$, $\varepsilon>0$. Hence our Standing Hypothesis
\ref{H} is satisfied, with $\beta_{\delta}(\varepsilon) =
Ce^{-\lambda_0 \delta/\varepsilon}$.

It remains to check Hypothesis \ref{hy.na}. The nonattractors in this
case are the unstable equilibria $(x^*_j)_{j=0,\ldots,k}$. Let
$\delta_1$ be the positive parameter associated with the stable points
$(s_i)_{i=1,\ldots,l}$ in our Corollary \ref{summary}. By their
assumption (A5), there exists $\beta>0$ such that we have, for
$j=1,\ldots,k$, $\varepsilon$ and $\delta$ small enough,
\[
\inf_{x \in V_j} \mathbb{P}_x \bigl(
\xi^{\varepsilon}(x) > \varepsilon \bigr) \geq\beta,\qquad \inf_{x \in V_j}
\mathbb{P}_x \bigl(\xi^{\varepsilon
}(x) < -\varepsilon \bigr) \geq
\beta,
\]
where $V_j = N^{\delta}(x^*_j)$. Call $\alpha_n = \delta/\varepsilon_n$
and assume without loss of generality that it is an integer. Using the
fact that $x^*_j$ is an unstable equilibrium, for $n$ large enough,
\[
\inf_{x \in V_j} \mathbb{P}_x \bigl(\bigl|X_{\alpha_n}^{\varepsilon
_n}-x^*_j\bigr|\geq\delta \bigr) \geq\beta^{\alpha_n}.
\]
Now choose $0<a<\lambda_0 \delta_1$. We have, by Markov's property,
\[
\sup_{x \in V_j} \mathbb{P}_x \bigl(
\tau_{V_j}^n \geq e^{a/\varepsilon
_n} \bigr) \leq\bigl(1-
\beta^{\alpha_n}\bigr)^{\nu_n},
\]
where $\nu_n = \frac{\varepsilon_n e^{a/\varepsilon_n}}{\delta}$. By
construction, $e^{a/\varepsilon_n} \beta_{\delta_1}(\varepsilon_n)$
goes to zero as $n$ goes to infinity. Additionally,
\[
\bigl(1-\beta^{\alpha_n}\bigr)^{\nu_n} \sim_{n \rightarrow+ \infty} \exp
\biggl(- \frac{\varepsilon_n e^{a/ \varepsilon_n} e^{\delta\log(\beta
)/\varepsilon_n}}{\delta} \biggr).
\]
This quantity vanishes as $n$ goes to infinity if we choose $\delta$
small enough (more precisely, $\delta$ must be chosen smaller than
$-a/\log\beta$). Therefore, we have verified Hypothesis \ref{hy.na},
and we can apply our result to conclude that the support of weak* limit
points of the QSDs do not include the unstable equilibria $x^*_j$.

%%%%%%%%%%%%%%%%%%%%%%%%%%%%%%%%%%%%%%%%%%%%%%%%%%%%%%%%%%%%%%%%%
%%%%%%%%%%%%%%%%%%%%%%%%%%%%%%%%%%%%%%%%%%%%%%%%%%%%%%%%%%%%%%%%%
%s4 #&#
%s4 ###
\section{\texorpdfstring{Properties of $\rho$-basic sets and proof of Theorem \protect\ref{th.ld}}
{Properties of rho-basic sets and proof of Theorem 2.7}}\label{s4}

In this section, we assume that the assumptions of Theorem~\ref{th.ld}
hold. We prove that Hypotheses~\ref{hy.rho},
\ref{hy.rho2}~and~finiteness of the $\rho$-basic classes imply
Hypotheses~\ref{hy.finite}, \ref{hy.na}~and~\ref{hy.boundary}. From
these implications, it follows that Theorem~\ref{th.suppmu},
Lemma~\ref{lm.lambda} and Proposition~\ref{prop.inv} imply that
Theorem~\ref{th.ld} holds. Indeed, by Lemma~\ref{lm.lambda}, we obtain
the lower bound for $\lambda_{\varepsilon}$ and, by
Proposition~\ref{prop.inv}, that any weak*-limit point of
$\{\mu_{\varepsilon}\}_{\varepsilon>0}$ is \mbox{$F$-}invariant. The
fact that these limiting measures sit on $\rho$-quasiattractors follows
from Theorem~\ref{th.suppmu}.

%
%le4.1 #&#
\begin{lemma}\label{lm.boundary}
Hypotheses \ref{hy.rho} and \ref{hy.rho2} imply Hypothesis
\ref{hy.boundary}.
\end{lemma}

\begin{pf} Let $\delta_0>0$ be given. Assertion (iii) of
Hypothesis~\ref{hy.rho} implies that
\[
c=\tfrac{1}{4} \inf\bigl\{\rho(x,y)\dvtx x,y\in M, d\bigl(F(x),y\bigr)>
\delta_0\bigr\}>0.
\]
It follows from the definition of $\beta_\delta$ and inequality
\eqref{upper} that $\beta_{\delta_0}(\varepsilon)\leq  \exp
(-3c/\varepsilon )$ for $\varepsilon>0$ sufficiently small.
Hypothesis~\ref{hy.rho2} implies that there exists an open neighborhood
$V_0$ of $M_0$ such that $p^\varepsilon(x,M_0)\geq
\exp(-2c/\varepsilon)$ for $\varepsilon>0$ sufficiently small and $x\in
V_0$. Hence,
\[
\lim_{\varepsilon\to0} \frac{\beta_{\delta_0}(\varepsilon)}{\inf_{x\in
V_0} p^\varepsilon(x,M_0)} \leq \lim_{\varepsilon\to0}
\exp (-c/\varepsilon) = 0.
\]\upqed
\end{pf}

For the remaining implications, we need to gain some insights about the
relation between ap and $\rho$-chain recurrence. As their ap
counterparts, the $\rho$-basic classes whose closure is in $M_1$ are
actually closed. This follows from the next two lemmas. We consider the
quantity
\[
\alpha^* = \sup_{\delta>0} \inf \bigl\{\rho(x,y) \dvtx x \in M, y \in
M, d\bigl(F(x),y\bigr) > \delta \bigr\} \in(0,+\infty].
\]

%
%le4.2 #&#
\begin{lemma} \label{lm.po=compact}
For any $0\leq\alpha< \alpha^*$, there exists $\delta>0$ such that, for
any $\xi=(\xi_0,\ldots,\xi_n)$ satisfying $A_n(\xi) < \alpha$, we have
$\|\xi_i\| \leq {\|F\|+\delta}$ for $i=1,\ldots,n$. In particular,
$\mathcal{R}_{\rho}$ is bounded.
\end{lemma}

\begin{pf} Given $0\leq\alpha< \alpha^*$, there exists $\delta>0$ such that
\[
\rho(x,y) < \alpha\Rightarrow d\bigl(F(x),y\bigr) < \delta.
\]
Therefore, if $A_n(\xi) < \alpha$, then $\rho(\xi_i,\xi_{i+1}) <
\alpha$ for $i=0,\ldots,n-1$, which implies that $d(\xi_{i+1},F(\xi_i))
<\delta$ for $i=0,\ldots,n-1$.
\end{pf}

%
%le4.3 #&#
\begin{lemma} \label{lm.B} We have the following:
\begin{longlist}[(iii)]
\item[(i)] The function $B_{\rho}$ is upper semicontinuous on $M_1
    \times M$.

\item[(ii)] Let $\eta\in(0, \alpha^*)$ and $y \in M$. If the set
    $\{x \in M \dvtx B_{\rho}(x,y) \leq\eta\}$ has its closure in
    $M_1$, then it is closed.

\item[(iii)] Let $x \in M_1$. Assume that the set $\{y \in M
    \dvtx B_{\rho}(x,y) \leq\eta\}$ has its closure in $M_1$ for
    $\eta$ small enough. Then there exists $\eta_0$ such that, for
    any $\eta< \eta_0$, $\{y \in M \dvtx B_{\rho}(x,y) \leq\eta\}$
    is closed.
\end{longlist}
\end{lemma}

\begin{pf} Part (i) is proved in \citeauthor{kif88} [(\citeyear{kif88}), Lemma~5.1, page~58]. It
relies on the continuity of $\rho$ on $M_1 \times M$.

For part (ii), let $\{x_k\}_{k\geq 1}$ be a sequence of points in
$M_1$, such that\break  $\lim_{k\to\infty}x_k = x \in M_1$ and
$B_{\rho}(x_k,y) \leq\eta$. Pick $r>0$ such that $\overline{N^r(x)}
\subset M_1$. For any $\gamma>0$ and any $k\geq 1$ there exists
$\xi^{\gamma,k}=(\xi_0^{\gamma,k},\ldots,\xi^{\gamma,k}_{n_k})$ such
that $\xi^{\gamma,k}_0 = x_k$, $\xi_{n_k}^{\gamma,k} = y$ and
$A_{n_k}(\xi^{\gamma,k}) \leq\eta+ \gamma$. By Lemma
\ref{lm.po=compact}, for $\gamma< \alpha^* - \eta$, there exists
$\delta>0$ such that $\|\xi^{\gamma,k}_1\| \leq \|F\|+\delta$. Since
$\overline{N^r(x)}$ is a compact set contained in $M_1$, and $\rho$ is
continuous on $M_1 \times M$, $\rho$ is uniformly continuous on
$\overline{N^r(x)} \times\overline{N^{\|F\|+\delta}(0)}$. Thus
\[
\lim_{k\to\infty} \bigl|\rho\bigl(x,\xi^{\gamma,k}_1
\bigr) - \rho\bigl(x_k, \xi^{\gamma,k}_1\bigr)\bigr|= 0.
\]
Therefore,
\begin{eqnarray*}
B_{\rho}(x,y)&\leq & \liminf_{k\to\infty} \bigl(\rho\bigl(x,
\xi_1^{\gamma,k}\bigr)+\rho\bigl(\xi_1^{\gamma,k},
\xi_2^{\gamma,k}\bigr)+\cdots+ \rho\bigl(\xi _{n_k-1}^{\gamma,k},y
\bigr) \bigr)
\\
&\leq & \liminf_{k\to\infty} \bigl|\rho\bigl(x,\xi_1^{\gamma,k}
\bigr)-\rho\bigl(x_k, \xi _1^{\gamma,k}\bigr)\bigr|+
A_{n_k}\bigl(\xi^{\gamma,k}\bigr) \leq \eta+\gamma.
\end{eqnarray*}
Since this holds for any $\gamma>0$, part (ii) follows.

Proof of part (iii) is similar. However, we have to be careful since
$\rho$ is not continuous on $M_0 \times M$. Let $x\in M_1$ be given,
and assume that there exists $\overline{\eta}>0$ such that $K =
\overline{\{y \in M \dvtx B_{\rho}(x,y) \leq\overline{\eta}\}} \subset
M_1$. Define $a = d(M_0 \cap N^{\|F\|}(0), K) >0$. Since $\rho$
is continuous on $M_1\times M$ and $\rho(z,y)=0$ if and only if
$y=F(z)$, there exists $\alpha>0$ such that
\[
\rho(z,y) < \alpha\Rightarrow d\bigl(F(z),y\bigr) < \min(a/2,1)
\]
for all $z,y\in M$. Let $\eta_0 = \min(\alpha,\overline{\eta})$ and
choose $0<\eta< \eta_0$. We claim that $\{y \in M \dvtx B_{\rho}(x,y)
\leq\eta\}$ is closed. To see why, let $\{y_k\}_{k\geq 1}$ be a
sequence in $M_1$ such that $\lim_{k\to\infty}y_k =y \in M_1$ and
$B_{\rho}(x,y_k) \leq\eta$ for all $k$. For any $\gamma>0$, there
exists a~family
$\xi^{\gamma,k}=(\xi_1^{\gamma,k},\ldots,\xi^{\gamma,k}_{n_k})$ such
that $\xi^{\gamma,k}_0 = x$, $\xi_{n_k}^{\gamma,k} = y_k$ and
$A_{n_k}(\xi^{\gamma,k}) \leq \eta+ \gamma$. For $\gamma< \eta_0
- \eta$, $\rho(\xi_{n_k-1}^{\gamma,k}, y_k)< \alpha$. Therefore,
$d(F(\xi_{n_k-1}^{\gamma,k}),\break  y_k) < \min(a/2,1)$, which implies
that $d(F(\xi_{n_k-1}^{\gamma,k}),M_0)>a/2$. By continuity of $F$ and
$F$-invariance of $M_1$, the sequence $\{\xi^{\gamma,k}_{n_k-1}\}_k$ is
bounded away from~$M_0$. Since $\eta_0 < \alpha$,\vadjust{\goodbreak} Lemma
\ref{lm.po=compact} implies that there exists $\delta>0$ such that
$\{\xi^{\gamma,k}_{n_k- 1}\}_k\subset N^{\|F\|+\delta}(0)$. The
remainder of the proof is as for part (ii), with
$\xi^{\gamma,k}_{n_k-1}$ playing the role of $\xi
_1^{\gamma,k}$.
\end{pf}

Since boundedness of $\rho$-basic classes follows from Lemma
\ref{lm.po=compact}, Lemma~\ref{lm.B} implies that given a $\rho$-chain
recurrent point $x$, if $\overline{[x]_{\rho}} \subset M_1$, then
$[x]_{\rho}$ is compact. Clearly, if $x$ is $\rho$-chain recurrent and
$[x]_{\rho}$ is closed, then $x$ is ap-chain recurrent and $[x]_{\rho}
\subset[x]_{\mathrm{ap}}$, but the converse is not true in general. For
example, consider a situation where the ap $\delta$-chains joining $x$
to itself have arbitrarily large length as $\delta$ goes to zero, in
which case we could have $x \sim_{\mathrm{ap}} x$ but $x \nsim_{\rho}
x$. While Remark \ref{rk.ap} holds for $\rho$-chain recurrence, it is
not immediate, and therefore we provide a short proof. In the sequel,
we will call $\delta$-$\rho$-pseudoorbit any family $\xi_0,\ldots,\xi_n$
such that $A_n(\xi) \leq \delta$.

%
%le4.4 #&#
\begin{lemma} \label{lm.po} Let $[x]_{\rho}$ be a closed $\rho$-basic
class in $M_1$. For any $\theta>0$, there exists $\delta>0$ such that
any $\delta$-$\rho$-pseudoorbit joining $[x]_{\rho}$ to itself is
contained in $N^{\theta}([x]_{\rho})$.
\end{lemma}

\begin{pf} Pick $\theta$ small enough so that
$\overline{N^{\theta}([x]_{\rho})} \subset M_1$. Since $F$ is
$M_1$-invariant and $\overline{N^{\theta}([x]_{\rho})}$ is compact and
contained in $M_1$, so is its image by $F$. Hence, by closedness of
$M_0$, there exists a compact set $K \subset M_1$, which contains the
$\gamma$-neighborhood of $F(N^{\theta}([x]_{\rho}))$, for some $\gamma
>0$. Assume, by contradiction, that there exist a decreasing sequence
$\delta_k \downarrow0$ and $\delta_k$-$\rho$-pseudoorbits
$(\xi^k_{0},\ldots,\xi^k_{n_k})$ [i.e., $A_{n_k}(\xi^k) \leq\delta_k$]
such that $ \lim_{k\to\infty}\xi^k_0=u \in[x]_{\rho}$,
$\lim_{k\to\infty}\xi^k_{n_k} = w \in[x]_{\rho}$ and $j_k = \min\{j
\geq1 \dvtx \xi^k_j \notin N^{\theta}([x]_{\rho}) \} < n_k$. For $k$
large enough, we have
\[
\rho(z,y) < \delta_k\quad\Rightarrow\quad d\bigl(F(z),y\bigr) < \gamma
\]
for all $z,y\in M$. Since $\xi^k_{j_k-1} \in N^{\theta}([x]_{\rho})$,
we have $\xi^k_{j_k} \in K$ for $k$ large enough. By passing to a
subsequence if necessary, $\xi^k_{j_k}$ converges to some point $v \in
K \setminus N^{\theta}([x]_{\rho})$. On the other hand, consider the
pseudoorbits $\tilde\xi^k=(\xi^k_0,\ldots, \xi^k_{j_k-1},v)$. They
satisfy
\[
\lim_{k\to\infty} A_{j_k} \bigl(\tilde\xi^k
\bigr) \leq  \lim_{k\to\infty} A_{n_k}\bigl(\xi^k
\bigr) + \lim_{k\to\infty} \bigl| \rho\bigl(\xi_{j_k-1}^k,
v\bigr)-\rho\bigl(\xi _{j_k-1}^k, \xi_{j_k}^k
\bigr)\bigr| =0
\]
due to $A_{n_k}(\xi^k)\leq \delta_k$ and by uniform continuity of
$\rho$ on $K\times K$. Hence, $x<_\rho v$. Similarly, one can show that
$v<_\rho x$. Consequently, $v \in[x]_\rho$, a contradiction.
\end{pf}

%
%le4.5 #&#
\begin{lemma} The $\rho$-basic classes $[x]_{\rho}$ closed in $M_1$ are
invariant:\break  $F([x]_{\rho}) = [x]_{\rho}$.
\end{lemma}

\begin{pf} \hspace*{-2.5pt}Since $F(x)>_{\mathrm{ap}}\hspace*{-0.5pt}x$, $F([x]_{\rho}) \subset[x]_{\rho}$ follows
if we prove that \mbox{$F(x) <_{\rho}\hspace*{-0.5pt}x$}. Since
$[x]_{\rho}$ is compact and contained in $M_1$ (see Lemma
\ref{lm.po=compact}), we can find a compact\vadjust{\goodbreak} set $K$ bounded away from
$M_0$, $\delta_k \downarrow0$ and a family of
$\delta_k$-$\rho$-pseudoorbits $\xi^k=(\xi_0^k,\ldots,\xi_{n_k}^k)$ in
$K$ such that $\xi_0^{k}=\xi^{k}_{n_k} = x$. Consider the family
$\tilde{\xi}^{k} = (F(x),\xi_2^{k},\ldots,\xi^{k}_{n_k -1},x)$. By
uniform\vspace*{1pt} continuity of $\rho$ on $K\times K$ and the fact that
$\rho(x,y)=0$ if and only if $y=F(x)$, we may assume by passing to a
subsequence if necessary that $\lim_{k\to\infty}\xi^k_1 = F(x)$. Hence
\[
\lim_{k\to\infty} A_{n_k-1}\bigl(\tilde\xi^k
\bigr) \leq \lim_{k\to\infty} A_{n_k}\bigl(\xi^k
\bigr)+\bigl|\rho\bigl(F(x),\xi_2^k\bigr)-\rho\bigl(
\xi_1^k,\xi_2^k\bigr)\bigr| = 0
\]
as $A_{n_k}(\xi^k)\leq \delta_k$ and using uniform continuity of $\rho$
on $K\times K$. Hence, \mbox{$F(x)<_\rho x$}.

For the inclusion $[x]_{\mathrm{ap}}\subset F([x]_{\mathrm{ap}})$, pick
$y \in [x]_{\rho}$ such that $y \neq F(y)$ (if~there is no such $y$,
there is nothing to prove). For $\delta_k \downarrow0$, choose a family
of \mbox{$\delta_k$-$\rho$-}pseudoorbits
$\xi^k=(\xi_0^k,\ldots,\xi_{n_k}^k)$ in $K$ such that
$\xi_0^{k}=\xi^{k}_{n_k} = y$. Passing to a~subsequence if necessary,
we can assume that $\lim_{k\to\infty} \xi^k_{n_k-1}=z \in K$. Clearly,
$F(z)=y$ and $z \sim_{\rho} y$. Hence, $[x]_\rho\subset F([x]_\rho)$.
\end{pf}

The following proposition is a straightforward consequence of
Proposition 5.1 in \citet{kif88}.
%
%pr4.6 #&#
\begin{proposition} \label{pr.rhoQA=A} Let $[x]_{\rho}$ be an isolated
$\rho$-quasiattractor in $M_1$. Then it is an attractor and $[x]_{\rho}
= \bigcap_{\eta>0} D_{\eta}$, where
\[
D_{\eta} = \bigl\{y \in M \dvtx B_{\rho}(x,y) < \eta\bigr\}.
\]
\end{proposition}

Define the maximum distance on $M^{n+1}$ by $d_n(\zeta,\xi) =
\max_{j=0,\ldots,n} d(\zeta_j,\xi_j)$ for $(\zeta_0,\ldots,\zeta_{n})$,
$(\xi_0,\ldots,\xi_{n})\in M^{n+1}$. The following theorem and lemma
are analogous to the statements of Theorem 5.2(a) and Lemma 5.3 in
\citeauthor{kif88} [(\citeyear{kif88}), pages 66~and~72, resp.]. We
provide a proof of Theorem~\ref{th.d_n} that slightly differs from the
proof of Kifer. The proof of Lemma~\ref{lm.KaN} follows directly from
Kifer's proof of his Lemma 5.3.

%
%th4.7 #&#
\begin{theorem} \label{th.d_n}Let $K \subset M_1$ be a compact set.
Given $\eta, \delta, N >0$, there exists $\varepsilon_0 >0$ such that
\[
\mathbb{P}_x \bigl[d_n\bigl(\bigl(X^{\varepsilon}_0,\ldots,X^{\varepsilon}_n\bigr),\xi\bigr) < \eta \bigr] \geq\exp
\biggl( - \frac{A_n(\xi) + \delta}{\varepsilon
} \biggr)
\]
for any $x \in K$, $\varepsilon< \varepsilon_0$, $n \leq N$ and $\xi =
(\xi_0,\ldots,\xi_n)\in K^{n+1}$ with $\xi_0 = x$.
\end{theorem}

\begin{pf} Analogously to Kifer's proof, we introduce the quantity
\[
n^K_{\gamma} = \sup \bigl\{\bigl|\rho(y,z) - \rho
\bigl(y',z'\bigr)\bigr| \dvtx y,y' \in K, d
\bigl(y,y'\bigr) \leq\gamma, d\bigl(z,z'\bigr) \leq
\gamma \bigr\}.
\]
By uniform continuity of $\rho$ on compact subsets of $M_1 \times M_1$,
\mbox{$\lim_{\gamma\to0}n^K_{\gamma}= 0$}. Let $\eta$, $\delta$, and
$N$ be given. Choose $0<\gamma< \eta$ such that $n^k_{\gamma} <
\delta/2N$ and \mbox{$N^\gamma(K)\subset M_1$}.\vadjust{\goodbreak} Now let $\xi= (\xi_0 =
x,\xi_1,\ldots,\xi_n)\in K^{n+1}$. By the uniform lower bound
(\ref{lower}) of Hypothesis~\ref{hy.rho} there exists a function
$g\dvtx]0,+ \infty[\mbox{ }\rightarrow\mbox{ }]0,+\infty[$ such that
\mbox{$\lim_{\varepsilon\to0} g(\varepsilon)=0$} and
\[
\varepsilon\log p^{\varepsilon}\bigl(x,N^\gamma(\xi_i)
\bigr) \geq-\inf_{y \in N^\gamma(\xi_i)} \rho(x,y) - g(\varepsilon)
\]
for any $x \in N^\gamma(K)$ and $1\leq  i\leq  n$.

Hence we have
\begin{eqnarray*}
&& \mathbb{P}_x \bigl[ d_n\bigl(X^{\varepsilon},\xi\bigr)< \eta\bigr]
\\
&&\qquad \geq \mathbb{P}_x \bigl[ d_n \bigl(X^{\varepsilon},\xi\bigr) <\gamma\bigr]
\\
&&\qquad = \int_{x_1 \in N^\gamma(\xi_1)} p^{\varepsilon}(x,dx_1) \cdots
\int_{x_n \in N^\gamma(\xi_n)} p^{\varepsilon}(x_{n-1},dx_n)
\\
&&\qquad \geq p^{\varepsilon}\bigl(x,N^\gamma(\xi_1)\bigr)
\prod_{i=1}^{n-1} \inf_{x_i \in N^\gamma(\xi_i)} p^{\varepsilon}
\bigl(x_i,N^\gamma(\xi_{i+1})\bigr)
\\
&&\qquad \geq \exp\Biggl[-\frac{1}{\varepsilon} \Biggl(\inf_{y\in N^\gamma(\xi_1)}\rho(x,y)
+ \sum_{i=1}^{n-1} \sup
_{x_i \in N^\gamma(\xi_i)} \inf_{y_i \in N^\gamma(\xi_{i+1})} \rho(x_i,y_i)
+ n g(\varepsilon) \Biggr) \Biggr]
\\
&&\qquad \geq \exp\biggl[-\frac{1}{\varepsilon} \bigl(A_n(\xi) + n g(
\varepsilon) + n n^K_{\gamma}\bigr) \biggr].
\end{eqnarray*}
The result follows by choosing $\varepsilon_0$ small enough so that
$g(\varepsilon) \leq\delta/2N$, for every $\varepsilon<
\varepsilon_0$.
\end{pf}

%
%le4.8 #&#
\begin{lemma} \label{lm.KaN} Let $K$ be a compact set in $M$ which does
not contain any entire semiorbits $\{F^i(x), i \in\mathbb{N} \}$. Then
there exists $a>0$ and $N \in\mathbb{N}$ (which depend on~$K$) such
that:
\begin{longlist}[(a)]
\item[(a)] for any sequence $\xi\in K^n$ with $n>N$, we have
    $A_n(\xi) > (n-N)a$;

\item[(b)] there exists $\varepsilon_0>0$ such that, for any $n
>N$ and any $0<\varepsilon< \varepsilon_0$,
\[
\sup_{x \in K} \mathbb{P}_x \bigl[
\tau^{\varepsilon}_{K} >n\bigr] \leq e^{-((n-N)a)/\varepsilon},
\]
where $\tau^{\varepsilon}_{K} = \inf \{j \geq0 \dvtx
X^{\varepsilon}_j \notin K  \}$.
\end{longlist}
\end{lemma}

Recall that $\omega(x)= \bigcap_{n\geq 1}
\overline{\bigcup_{m\geq n} F^m(x)}$ and that a point $x \in M$ is
called \textit{nonwandering} if for all open neighborhoods $U$ of $x$
and any $N \in\mathbb{N}$, there exists $n \geq N$ such that $F^n(U)
\cap U \neq\varnothing$. We denote by $\mathrm{NW}(F)$ the set of
nonwandering points of $F$. Note that $\omega$-limit points are always
nonwandering: $ \{y \in M\dvtx y \in\omega(x), \mbox{for some } x  \}
\subset \mathrm{NW}(F)$.

%
%le4.9 #&#
\begin{lemma} \label{lm.omega} The set $\mathrm{NW}(F) \cap M_1$ is contained in
$\mathcal{R}_{\rho}$. In particular, any $\omega$-limit point in $M_1$
is also in $\mathcal{R}_{\rho}$.
\end{lemma}

\begin{pf} Let $x \in M_1 \cap \mathrm{NW}(F)$ and $\delta>0$ be given. By continuity
of $\rho$ in \mbox{$M_1 \times M$} and $\rho(x,F(x))=0$ for all $x$,
there exists $\gamma>0$ such that
\[
\rho\bigl(x,F(y)\bigr) < \delta/2 \quad\mbox{and}\quad \rho(z,x) < \delta/2
\]
for $y \in N^{\gamma}(x)$ and $F(z) \in N^{\gamma}(x)$. Since $x$ is
nonwandering, there exists $n\geq 1$ such that
\[
F^{n}\bigl(N^{\gamma}(x)\bigr) \cap N^{\gamma}(x) \neq
\varnothing.
\]
Pick $y, z \in N^{\gamma}(x)$ such that $F^{n}(y) = z$. Now consider
the chain $\xi= (x,F(y),\ldots,\break  F^{n-1}(y),x)$. Since
$F(F^{n-1}(y)) = z \in N^\gamma(x)$, we have $A(\xi) = \rho(x,F(y)) +
\rho(F^{n-1}(y),x) < \delta$. Taking $\delta\downarrow0$ yields
$x\sim_\rho x$ as claimed.
\end{pf}

%
%co4.10 #&#
\begin{corollary} \label{co.support} Assume that $\mu$ is an
$F$-invariant probability measure whose support $S$ lies in $M_1$. Then
$S \subset\mathcal{R}_{\rho}$.
\end{corollary}

\begin{pf} By the Poincar\'{e} recurrence theorem, $S$ is included in
the set
\[
\overline{\bigl\{x \in M_1 \dvtx x \in\omega(x)\bigr\}} \subset \mathrm{NW}(F).
\]
Applying Lemma~\ref{lm.omega} completes the result.
\end{pf}

%
%co4.11 #&#
\begin{corollary} \label{co.partsol} Assume that $\mathcal{R}_{\rho}
\cap M_1$ admits a neighborhood $U$, whose closure lies in some compact
set $K \subset M_1$. Then there exists $N \in\mathbb{N}$ such that any
partial solution $\zeta= (x,F(x),\ldots,F^n(x))\in K^{n+1}$ with $n
\geq N$ must pass through~$U$.
\end{corollary}

\begin{pf} The set $K \setminus U$ does not contain any entire
semiorbit of $F$, by Lemma~\ref{lm.omega}. Since $A_n(\zeta)=0$,
applying Lemma~\ref{lm.KaN}(a) completes the proof.
\end{pf}

We already stated that if $[x]_{\rho}$ is a closed $\rho$-basic class,
then $[x]_{\rho} \subset[x]_{\mathrm{ap}}$. Under the finiteness
assumption of Theorem \ref{th.ld}, we derive the equality between ap
and $\rho$-basic classes. Call $K_1,\ldots,K_v$ the $\rho$-basic
classes in $M_1$ (recall that they are supposed to be closed), and
label $K_1,\ldots,K_\ell$ the quasi-attractors among them. Proposition
\ref{pr.rhoQA=A} implies that $K_1,\ldots,K_\ell$ are attractors. The
following lemma implies that finiteness of the $\rho$-basic classes in
$M_1$ implies finiteness of the $ap$-basic classes in $M_1$. In
particular, Hypothesis \ref{hy.finite} holds under the assumptions of
Theorem~\ref{th.ld}.

%
%th4.12 #&#
\begin{theorem}\label{th.equalitybasic} Assume that
there is a finite number of $\rho$-basic classes in~$M_1$. Then
$\mathcal{R}_{\rho} \cap M_1=\mathcal{R}_{\mathrm{ap}}\cap M_1$ and
$[x]_{\rho} = [x]_{\mathrm{ap}}$ for any $x
\in\mathcal{R}_{\mathrm{ap}}\cap M_1$.
\end{theorem}

\begin{pf} Let $x \in\mathcal{R}_{\mathrm{ap}} \cap M_1$. We prove that $x \in
\mathcal{R}_{\rho} \cap M_1$ and $[x]_{\mathrm{ap}} \subset[x]_{\rho}$.
If $[x]_{\mathrm{ap}} = \{x\}$, then $x$ is a fixed point and there is
nothing to left to prove. Let $y \in[x]_{\mathrm{ap}}$, $y \neq x$ and
$\alpha>0$.\vadjust{\goodbreak}

Remark~\ref{rk.ap} implies there exists a compact set $K \subset M_1$
with $\bigcup_i K_i \subset K$, a~sequence $\delta_k \downarrow_k 0$
and a family $\xi^k=(\xi^k_0=x,\ldots,\xi_{n_k}^k=y)\in K^{n_k+1}$ of
ap \mbox{$\delta_k$-}pseudoorbits joining $x$ to $y$.

Let $\gamma>0$ be chosen so that $B_\rho(a,b) < \alpha$ for all
$i=1,\ldots,v$ and $a,b \in N^{\gamma}(K_i)$, the closure of $U =
\bigcup_i N^{\gamma}(K_i)$ is contained in $K$, and $N^\gamma(K_i)$ is
an isolating neighborhood for $K_i$ for all $i$. Corollary
\ref{co.partsol} implies that there exists a positive integer $N$ such
that every partial solution $\{a,F(a),\ldots,F^n(a)\}$ in $K$ of length
$n\geq  N$ must pass through $U$. Consequently, by compactness of $K^N$
and continuity of $F$, we can find $k_0$ such that $\xi^k$ cannot have
more than $N$ consecutive terms in $K \setminus U$ for $k\geq  k_0$.

Now, given $k \in\mathbb{N}$, define $\sigma_0(k) = 0$ and $\tau_0(k) =
\min\{j>0 \dvtx \xi_j^k \notin U\}$. For $i\geq 1$, define inductively
the terms $\sigma_i(k)=\min\{ j> \tau_{i-1}(k)\dvtx \xi^k_j \in U\}$
and $\tau_i(k)=\min\{j>\sigma_{i}(k) \dvtx \xi_j^k \notin U\}$. This
defines two sequences $\{\tau_i(k)\}_{i=0,\ldots,p_k}$ and
$\{\sigma_i(k)\}_{i=0,\ldots,q_k}$. Notice that $q_k = p_k$ if $y
\notin \bigcup_i K_i$ and $q_k = p_k +1$ otherwise. By truncating
multiple entries of ap pseudoorbits into each set $N^\gamma(K_i)$, we
can assume that $q_k \leq v-1$. After truncation, these pseudoorbits
may only satisfy $d(F(\xi_j^k),\xi_{j+1}^k) \leq\delta_k$ for
$\tau_i(k)-1\leq j \leq \sigma_{i+1}(k)$. Therefore
\begin{eqnarray*}
B_\rho(x,y) &\leq&\sum
_{i=0}^{p_k} \bigl(B_\rho\bigl(
\xi^k_{\sigma_i(k)},\xi ^k_{\tau_i(k)-1}\bigr) +
B_\rho\bigl(\xi^k_{\tau_i(k)-1},\xi^k_{\tau_i(k)}
\bigr)
\\
&&\hspace*{16pt}{} +B_\rho\bigl(\xi^k_{\tau_i(k)},
\xi^k_{\sigma_{i+1}(k)-1}\bigr) + B_\rho\bigl(
\xi^k_{\sigma_{i+1}(k)-1},\xi^k_{\sigma_{i+1}(k)}\bigr)\bigr)
\\
&&{}+ B_\rho\bigl(\xi^k_{\tau_{q_k}(k)},y\bigr)
\end{eqnarray*}
in the case where $y \in\bigcup_i K_i$, and
\begin{eqnarray*}
B_\rho(x,y) &\leq&\sum
_{i=0}^{p_k} \bigl( B_\rho\bigl(
\xi^k_{\sigma_i(k)},\xi ^k_{\tau_i(k)-1}\bigr) +
B_\rho\bigl(\xi^k_{\tau_i(k)-1},\xi^k_{\tau_i(k)}
\bigr) \bigr)
\\
&&{}+ \sum_{i=0}^{p_k-1} \bigl(B_\rho
\bigl(\xi^k_{\tau_i(k)}, \xi^k_{\sigma
_{i+1}(k)-1}\bigr)
+B_\rho\bigl(\xi^k_{\sigma_{i+1}(k)-1},\xi^k_{\sigma
_{i+1}(k)}\bigr) \bigr)
\\
&&{} + B_\rho\bigl(\xi^k_{\sigma_{p_k}(k)},y\bigr)
\end{eqnarray*}
otherwise. In either case, our choice of $\gamma$ implies
\[
B_\rho(x,y) \leq v \alpha+ \bigl(v (N+2)+1\bigr) \sup\bigl\{
\rho(a,b) \dvtx d\bigl(F(a),b\bigr) \leq\delta_k, a,b\in K \bigr\}
\]
for $k$ sufficiently large. Uniform continuity of $\rho$ on $K \times
K$ implies that $\lim_{k\to\infty} \sup\{ \rho(a,b) \dvtx
d(F(a),b) \leq\delta_k, a,b\in K \}=0$, and we obtain that\break
$B_\rho(x,y) \leq v \alpha$. Since this holds for any $\alpha>0$ we get
that $x <_{\rho} y$. Similarly, $y <_{\rho} x$, which yields $x
\sim_{\rho} y$. Therefore, $x \in\mathcal{R}_{\rho}$ and $[x]_{\rho} =
[x]_{\mathrm{ap}}$.
\end{pf}

The next proposition shows that Hypothesis \ref{hy.na} holds.

%
%pr4.13 #&#
\begin{proposition} \label{pr.NQA} Let $j \in\{\ell+1,\ldots,v\}$. We
can find $\eta>0$ such that, for any $\gamma>0$, there exists
$\varepsilon_0 >0$ (which depends on $\eta$ and $\gamma$) and a
function $\zeta$ on $(0,\varepsilon_0)$ such that $ \lim_{\varepsilon
\to0} \zeta(\varepsilon) = 0 $ and
\[
\sup_{x \in N^{\eta}(K_j)} \mathbb{P}_x\bigl[
\tau^{\varepsilon}_{N^{\eta
}(K_j)} > e^{\gamma/\varepsilon}\bigr] \leq\zeta(
\varepsilon)
\]
for any $\varepsilon< \varepsilon_0$.
\end{proposition}

\begin{pf} First of all, by definition of a non$\rho$-quasiattractor, there
exists $\eta>0$ such that the closure of $N^{2\eta}(K_j)$ belongs to
$M_1$, and for any $\gamma>0$ and any $x \in N^{\eta}(K_j)$, there
exists a sequence $\xi^{\gamma} =
(\xi^{\gamma}_0,\ldots,\xi^{\gamma}_{n(\gamma)})$ such that
\[
\xi^{\gamma} = x,\qquad \xi^{\gamma}_{n(\gamma)} \notin
N^{2\eta}(K_j) \quad\mbox{and}\quad A_{n(\gamma)}\bigl(
\xi^{\gamma}\bigr) < \gamma.
\]
Call $U = N^{2\eta}(K_j)$. Since $M_1$ is invariant by $F$,
$F(\overline{U})$ is compact and contained in $M_1$. Hence there
exists $r>0$ and a compact set $K \subset M_1$ such that
\[
N^{r}\bigl(F(\overline{U})\bigr) \subset K.
\]
By continuity of $\rho$ on $\overline{U} \times M$ and since $\rho$ is
strictly positive on the compact set $\overline{U} \times
(N^r(F(\overline{U})))^c$, there exists $\gamma_0>0$ such that
\[
\rho(x,y) > \gamma_0\qquad\mbox{for all } x \in\overline{U}, y \in K^c.
\]
In particular, this means that, for $\gamma< \gamma_0$, the sequence
$\xi^{\gamma}$ must pass through $K \setminus U$, and we can therefore
assume without loss of generality that $\xi^{\gamma}_{n(\gamma)} \in K
\setminus U$ and $\xi^{\gamma}$ lives in $K$.

Pick $\delta>0$. We now apply Theorem \ref{th.d_n} in the compact set
$K$, with $\delta$, $\eta$ and $N = n(\gamma)$: there exists
$\varepsilon_0>0$ [which depends on $\eta$, $\delta$ and $n(\gamma)$]
such that, for any $\varepsilon< \varepsilon_0$,
\[
\mathbb{P}_x \bigl[d_{n(\gamma)}\bigl(X^{\varepsilon},
\xi^{\gamma}\bigr) < \eta \bigr] \geq\exp \biggl( - \frac{\gamma+ \delta}{\varepsilon}
\biggr).
\]
Consequently there exists $\varepsilon_0^{\prime}>0$ such that, for any
$0<\varepsilon< \varepsilon_0^{\prime}$ [up to changing slightly
$n(\gamma)$], we have
\[
\mathbb{P}_x \bigl[ \tau^{\varepsilon}_{N^{\eta}(K_j)} \leq n(
\gamma) \bigr] > e^{-\gamma/\varepsilon}.
\]
Consequently,
\[
\mathbb{P}_x \bigl[ \tau^{\varepsilon}_{N^{\eta}(K_j)} \geq
e^{2\gamma
/\varepsilon} \bigr] \leq\bigl(1 - e^{-\gamma/\varepsilon}\bigr)^{ [
e^{2\gamma/\varepsilon}/n(\gamma) ]}.
\]
The last quantity is of order $\exp
(-e^{\gamma/\varepsilon}/(2n(\gamma))  )$ and therefore goes
to zero.
\end{pf}

By assumption of Theorem \ref{th.ld}, there is at least one
$\rho$-quasiattractor (which turns out to be an attractor by
Proposition \ref{pr.rhoQA=A}) among the $\rho$-basic classes included
in $M_1$, such that $\mu_n(U) >0$ for all $n$, for all open
neighborhoods $U$. Hence, by Lemma \ref{lm.lambda} there exists
$\delta_0>0$ such that $\lambda_n > 1 - \beta_{\delta_0}(\varepsilon_n)
> 1 - e^{-c_0/\varepsilon_n}$, where
\[
c_0 = \tfrac{1}{2} \inf \bigl\{\rho(x,y) \dvtx (x,y) \in
M_1 \times M, d\bigl(F(x),y\bigr) \geq\delta_0 \bigr\}
>0.
\]
Let $\delta_1 < \delta_0$ and $\{V_i\}_{i=1,\ldots,v}$ be chosen so
that Corollary \ref{summary} holds.

We are now ready to prove Theorem \ref{th.ld}. Since Hypotheses
\ref{hy.finite} and \ref{hy.boundary} are satisfied, it remains to
verify Hypothesis~\ref{hy.na}. Choose the neighborhoods $\{V_i\}_{
i=\ell+1,\ldots,v}$ such that $V_i \subset N^{\eta}(K_i)$, where $\eta$
is given by Proposition \ref{pr.NQA}. Choose $\gamma< \frac{c_1}{2}$
where
\[
0 < c_1 = \inf \bigl\{\rho(x,y)\dvtx (x,y) \in M_1 \times
M, d\bigl(F(x),y\bigr) \geq \delta_1 \bigr\} \leq 2 c_0.
\]

Define $h(\varepsilon) = e^{\gamma/\varepsilon}$.
Proposition~\ref{pr.NQA} implies that there exists $\varepsilon_0>0$
and a function $\zeta$ such that $\lim_{\varepsilon\to0}
\zeta(\varepsilon)=0$ and
\[
\sup_{x \in N^{\eta}(K_j)} \mathbb{P}_x \bigl[
\tau^{\varepsilon_n}_{N^{\eta
}(K_j)} > h(\varepsilon_n)\bigr] \leq
\zeta(\varepsilon_n) \rightarrow0
\]
for any $\varepsilon< \varepsilon_0$. Since $\lim_{n\to\infty}
h(\varepsilon_n) \beta_{\delta_1}(\varepsilon_n) = 0$, Hypothesis
\ref{hy.na} holds. This completes the proof of Theorem~\ref{th.ld}.

%%%%%%%%%%%%%%%%%%%%%%%%%%%%%%%%%%%%%%%%%%%%%%%%%%%%%%%%%%%%%%%
%%%%%%%%%%%%%%%%%%%%%%%%%%%%%%%%%%%%%%%%%%%%%%%%%%%%%%%%%%%%%%%
%s5 #&#
%s5 ###
\section{Applications}\label{s5}

Our results are broadly applicable to many Markov chain models in
population biology. To give some flavor of this applicability, we
introduce two classes of Markov chains satisfying our probabilistic
assumptions and some illustrative applications to metapopulation
dynamics, competing species, host-parasitoid interactions and
evolutionary games. For each application there are two ingredients for
verifying the conditions of Theorem~\ref{th.ld}. The probabilistic
ingredient involves verifying that there exist quasi-stationary
distributions and verifying the large deviation assumptions. We defer
verifying these conditions until Section~\ref{s6}. The dynamical
ingredient involves verifying there is a finite number of $\rho$-basic
classes and identifying the attractors. For the second ingredient, we
introduce a proposition, that is, applicable to most of our examples.

%s5.1 #&#
%s5.1 ###
\subsection{A dynamical proposition}\label{s5.1}
To state the proposition, we need a few definitions from dynamical
systems. For $x\in M$, let $\omega(x)=\{y$: there exists $n_k\to\infty$
such that $\lim_{k\to\infty}F^{n_k}(x)=y\}$ be the
\textit{$\omega$-limit set for $x$} and $\alpha(x)=\{y$: there exist
$n_k\to\infty$ and $y_k \in M$ such that $F^{n_k}(y_k)=x$ and
$\lim_{k\to\infty} y_k =y\}$ be the \textit{$\alpha$-limit set for
$x$}. Our assumption that $F$ is bounded implies that there exists a
global attractor given by the compact, $F$-invariant set
$\Lambda=\bigcap_{n\geq 0} F^n(M)$. For all $x\in\Lambda $, $\omega(x)$
and $\alpha(x)$ are compact, nonempty, $F$-invariant sets.

A \textit{Morse decomposition} of the dynamics of $F$ is a collection
of $F$-invariant, compact sets $K_1,\ldots,K_k$ such that:
\begin{itemize}
\item $K_i$ is isolated, that is, there exists a neighborhood of
    $K_i$ such that it is the maximal $F$-invariant set in the
    neighborhood, and

\item for every $x\in\Lambda\setminus\bigcup K_i$, there exist $i>j$
such that $\omega(x) \subset K_i$ and \mbox{$\alpha(x)\subset K_j$}.
\end{itemize}
Modulo replacing the invariant sets $K_i$ by points, one can think of
$F$ being gradient-like as all orbits move from lower indexed invariant
sets to higher indexed invariant sets.

%
%pr5.1 #&#
\begin{proposition}\label{prop.dyn} Assume Hypothesis~\ref{hy.rho} holds.
If $F$ admits a Morse decomposition $K_1,\ldots,K_k$ such that:
\begin{itemize}
\item $K_i \subset M_1$ or $K_i\subset M_0$ for each $i$, and

\item $K_i$ is transitive whenever $K_i\subset M_1$, that is, there
    exists $x\in K_i$ such that $\{x,F(x),F^2(x),\ldots\}$ is dense
    in $K_i$,
\end{itemize}
then $\rho$-basic classes in $M_1$ are given by the $K_i \subset M_1$.
In particular, there is a finite number of $\rho$-basic classes in
$M_1$, and each of them is closed.
\end{proposition}

\begin{pf} $\!\!\!$Let $K_1,\ldots,K_k$ be a Morse decomposition for $F$. Let
$I\subset\{1,\ldots,k\}$ be such that $K_i\subset M_1$ if and only
$i\in I$. By assumption, $K_i\subset M_0$ for $i\notin I$, and $K_i$ is
transitive for $i\in I$. Transitivity of $K_i$ for $i\in I$ and
continuity of $\rho$ restricted to $M_1\times M$ implies that $K_i$ is
contained in a $\rho$-basic class for $i\in I$, that is, $x\sim_\rho y$
for all $x,y\in K_i$. As shown in Section~\ref{s4}, assertion (iii) of
Hypothesis~\ref{hy.rho} implies that the $\rho$-basic classes are
contained in the ap-chain recurrent set of $F$ which is contained in
$\bigcup_i K_i$. Hence, the $\rho$-basic classes in $M_1$ are given by
$\{K_i\}_{i\in I}$.
\end{pf}

%%%%%%%%%%%%%%%%%%%%%%%%%%%%%%%%%%%%%%%%%%%%%%%%%%%%%%%%%%%%%%%%%%%%%%%%%%%%%%%%%%%%%%%%%%%%%%
%s5.2 #&#
%s5.2 ###
\subsection{Nonlinear Poisson branching processes}\label{s5.2}

To describe structured populations with $k$ types of individuals (e.g.,
different genotypes or species, individuals living in different spatial
locations), let $x$ represent the vector of population densities
which lies in the nonnegative cone $\R^k_+$ of $\R^k$. A widely used
class of models in population biology~\citep{caswell-01} is the
nonlinear matrix model of the form $F(x)=A(x)x$ where $A(x)$ are
nonnegative matrices representing transitions births, deaths and
transitions between types of individuals (e.g., due to mutation or
dispersal).

Since real populations involve finite numbers of individuals, these
deterministic models can be viewed as approximations of more realistic,
stochastic representations of the population dynamics. More
specifically, let $1/\varepsilon$ that represents the ``size'' (e.g.,
area, volume) of the habitat. Let $N_t^\varepsilon\in\mathbb{Z}^k_+$
denote the vector of population abundances where $\mathbb{Z}^k_+$ is
the nonnegative cone of the $k$-dimensional integer lattice. Then
$X_t^\varepsilon=\varepsilon N_t^\varepsilon\in\varepsilon
\mathbb{Z}^k_+$ is the vector of population densities. For every $x\in
\R^k_+$, let $Z_1(x),Z_2(x),\ldots$ be a sequence of i.i.d. random
vectors with independent components, and whose $i$ component has a
Poisson distribution with mean~$x_i$. Given $N_0^\varepsilon\in
\mathbb{Z}^k_+$, we can define the Markov chains $\{X_t^\varepsilon\}$
iteratively by
\[
N_{t+1}^\varepsilon= Z_{t+1}\bigl(A\bigl(X_t^\varepsilon
\bigr) N_{t}^\varepsilon\bigr)\quad\mbox{and}\quad X_{t+1}^{\varepsilon}=
\varepsilon N_{t+1}^{\varepsilon}.
\]
Equivalently, we can write $X_{t+1}^\varepsilon= \varepsilon
Z_{t+1}(F(X_t^\varepsilon)/\varepsilon)$.

A useful observation about these Poisson processes, from the modeling
standpoint, is that multinomial sampling of a Poisson process still
corresponds to a Poisson process. More specifically, consider a
multinomial random vector $(X_1,\ldots, X_k)$ where the number of
samples $N$ is Poisson distributed with mean $\lambda>0$ and the
sampling probabilities are $(p_1,\ldots, p_k)$. Then
\begin{eqnarray*}
&& \mathbb{P}[X_1=x_1,\ldots,X_k=x_k]
\\
&&\qquad  =\mathbb{P}[X_1=x_1,\ldots,X_k=x_k \mid N=x_1+\cdots +x_k]P[N=x_1+\cdots+x_k]
\\
&&\qquad =\frac{(x_1+\cdots+x_k)!}{x_1!\cdots x_k!} p_1^{x_1} \cdots
p_k^{x_k} \exp(-\lambda) \frac{\lambda^{x_1+\cdots+x_k}}{(x_1+\cdots+ x_k)!}
\\
&&\qquad = \prod_{i=1}^k \frac{ (p_i \lambda)^{x_i}}{x_i!}\exp(-p_i \lambda).
\end{eqnarray*}
Hence, $X_1,\ldots,X_k$ are independent Poisson random variables with
rate parameters $p_1 \lambda,\ldots, p_k \lambda$. We make repeated use
of this observation in the examples provided below.

Since $F$ is bounded, one can show quite generally that these Markov
chains support quasi-stationary distributions whenever there are
absorbing sets. For all of our examples, these absorbing sets are
$\{0\}$ or $\partial\R^k_+=\{x\in\R^k_+\dvtx\prod x_i =0\}$. A~proof of
this assertion is given in Proposition \ref{prop.quasi} of
Section~\ref{s6}. Under slightly stronger assumptions (namely $A$ is
continuous and $F_i$ is strictly positive), we show in Proposition
\ref{pr.largedev} of Section~\ref{s6} that these Poisson processes also
satisfy our large deviation Hypotheses~\ref{hy.rho} and \ref{hy.rho2}.

To provide a taste of the possible applications, we apply our results
to three particular classes of nonlinear Poisson branching processes.

\subsubsection*{Metapopulation dynamics}
A fundamental question in population biology is how do local
demographic processes, such as reproduction and survivorship, interact
with dispersal (a regional demographic process) to determine
spatial-temporal patterns of abundance~\citep
{earn-etal-00,hastings-botsford-06,earn-levin-06,prsb-10}. This issue
has been studied extensively with discrete-time deterministic models
representing space as a finite collection of patches connected by
dispersal. To illustrate how our results apply to these metapopulation
models, we introduce a stochastic version of the spatial Ricker map,
which was originally studied by \citet{hastings-93} for $2$
patches, and for which we allow an arbitrary number, $k$, of patches.

Let $1/\varepsilon>0$ be the area or volume of a single patch,
$N_t^{\varepsilon,i}$ denote the number of individuals in patch $i$,
$N_t^\varepsilon=(N_t^{\varepsilon,1},\ldots,N_t^{\varepsilon,k})$ the
vector of population abundances across space and $X_t^\varepsilon
=\varepsilon N_t^\varepsilon$ the vector of population densities. To
describe reproduction within a patch, let $f(x)=f_0 \exp(-x)$ be the
mean fecundity of an individual when the local population density is
$x$ and the ``intrinsic'' fitness is $f_0>0$. The map $x\mapsto x f(x)$
is known as the Ricker map in theoretical ecology and is commonly used
to describe the population dynamics of a single
species~\citep{ricker-54,hastings-97,wysham-hastings-08}. Let
$D=(d_{ij})$ be an irreducible, row-stochastic matrix where $d_{ij}$
corresponds to the probability of an individual dispersing from patch
$i$ to patch $j$. Given $N_t^\varepsilon$, we define the
\textit{spatial Ricker process} as follows:
\begin{itemize}
\item Each individual in patch $i$ independently produces a
    Poisson-distributed number of offspring with mean
    $f(X_t^{\varepsilon,i})$ to replace themselves. Let
    $Z_{t+1}^{i,\varepsilon}$ be the total number of offspring
    produced in patch $i$, which is Poisson distributed with mean
    $N_t^{\varepsilon,i}f(X^{\varepsilon,i}_t)$. We assume that the
    $Z_{t+1}^{1,\varepsilon},\ldots,Z_{t+1}^{k,\varepsilon}$ are
    independent; there are no correlations in the reproductive
    output between distinct patches.

\item Independent of one another, offspring in patch $i$ move to
    patch $j$ with probability $d_{ij}$. To represent this
    movement, let
    $W_{t+1}^{\varepsilon}(i)=(W_{t+1}^{\varepsilon,1}(i),\ldots,\break W_{t+1}^{\varepsilon,k}(i))$
    be a multinomial random vector with sampling probabilities
    $d_{i1},\ldots,d_{ik}$ and $Z_{t+1}^{\varepsilon,i}$ trials.

\item Define
\[
N_{t+1}^\varepsilon= \sum_iW_{t+1}^{\varepsilon}(i)
\quad\mbox{and}\quad X_{t+1}^\varepsilon = \varepsilon N_{t+1}^\varepsilon.
\]
\end{itemize}
By our earlier observation about multinomial sampling of a Poisson
random variable, $\{X^\varepsilon_t\}$ is a nonlinear Poisson process
with
\[
F_i(x)=\sum_j d_{ji}
x_j f(x_j).
\]
Since $\| F \| \leq  f_0 \|D\|$ and $D$ is irreducible,
Proposition~\ref{prop.quasi} implies the process $\{X_t^\varepsilon\}$
has a quasi-stationary distribution $\mu_\varepsilon$ with respect to
the absorbing state $M_0=\{0\}$ for all $\varepsilon>0$.

Let $\mu$ be a weak* limit point of $\mu_\varepsilon$ as $\varepsilon
\to0$. To say something about the support of $\mu$, we need to
understand the dynamics of the map $F(x)$. The simplest applicable
result is a persistence and extinction dichotomy. Since $F_i(x)\leq f_0
x_i$, it follows that $0$ is a global attractor for $F(x)$ whenever
$f_0<1$. Alternatively, when $f_0>1$, a result of
\citeauthor{kon-etal-04} [(\citeyear{kon-etal-04}), Theorem~3] implies
that $F(x)$ has a positive attractor. Lemma~\ref{lm.lambda}(b$'$) and
Theorem~\ref{th.M0} imply the following result.
%
%pr5.2 #&#
\begin{proposition} Let $\mu$ be a weak* limit point for
quasi-stationary distributions $\mu_\varepsilon$ of the spatial Ricker
process $\{X_t^\varepsilon\}$. Then:
\begin{longlist}
\item[Extinction.] If $f_0<1$, then $\mu(\{0\})=1$.

\item[Metastability.] If $f_0>1$, then there exists a $\delta>0$
    such that \mbox{$\mu(N^\delta(\{0\}))=0$}.
\end{longlist}
\end{proposition}

In the limiting cases where the population is either weakly mixed or
well mixed, we can say more about the support of the limiting measure
$\mu$. These stronger assertions rely on the one-dimensional map
$x\mapsto f_0 x \exp(-x)$ having a linearly stable periodic orbit, call
it $\mathcal{S}=\{p,F(p),\ldots, F^{n-1}(p)\}$ where $n$ is the period.
\citeauthor{kozlovski-03} [(\citeyear{kozlovski-03}), Theorem~C] proved
that, for an open and dense set of $f_0$ values, such a stable periodic
orbit exists. Hence, this assumption is not very restrictive.

%
%th5.3 #&#
\begin{theorem}\label{th5.3}
Assume the one-dimensional map $x\mapsto f_0 x \exp(-x)$ has a linearly
stable periodic orbit, call it $\mathcal{S}=\{p,F(p),\ldots,
F^{n-1}(p)\}$, and $D$ is an irreducible, nonnegative matrix whose row
sums equal one (i.e., a row stochastic matrix).
\begin{longlist}
\item[Weakly mixing.] If $D$ is sufficiently close to the identity
matrix, then there exists $n_k$ linearly stable periodic orbits
for $F(x)$, and $\mu$ is supported by the union of these stable
periodic orbits.

\item[Strongly mixing.] If all the entries of $D$ are sufficiently
close to $1/k$ and the column sums of $D$ equal one (i.e., $D$
is doubly stochastic), then there exists a unique globally
stable periodic orbit for $F(x)$, and the support of $\mu$ is
given by this periodic orbit.
\end{longlist}
\end{theorem}

We remark that in the special case of a single patch, $k=1$, we recover
results of \citet{hognas-97},
\citet{klebaner-lazar-zeitouni-98} and
\citet{ramanan-zeitouni-99} for one-dimensional maps on a compact
interval. See Section~\ref{s3.5} for further discussion about this
point.

\begin{pf*}{Proof of Theorem~\ref{th5.3}} To prove the first assertion, consider
the uncoupled map
\[
\widetilde F(x)=\bigl(x_1f(x_1),x_2f(x_2),\ldots,x_kf(x_k)\bigr).
\]
Each of the components of this limiting map are given by the
one-dimensional map $g(x_i)=x_if(x_i)$ which by assumption has a
linearly stable periodic orbit $\mathcal{S}=\{p,g(p),\ldots,
g^{n-1}(p)\}$. This linearly stable periodic orbit gives rise to $n^k$
periodic orbits of the form $(g^{n_1}(p),\ldots,g^{n_k}(p))$ with
$0\leq n_j<n$ for $\widetilde F$. Since $g$ has a negative Schwartzian
derivative and a single critical point, \citeauthor{vanstrien-81}
[(\citeyear{vanstrien-81}), Theorem~A] proved that the complement of
the basin of attraction of ${\mathcal S}$ for $g$ can be decomposed
into a finite number of compact, $g$-invariant sets which have a dense
orbit and are hyperbolic repellers: there exists $C>0$ and $\lambda>1$
such that $|(g^n)'(x)|\geq C \lambda^n$ for all points $x$ in the set
and $n\geq 1$. Consequently, the $k$-dimensional mapping $\tilde F$ is
an Axiom A endomorphismn [\citet{przytycki-76}, page~271]: the derivative of $\widetilde
F$ is nonsingular for all points in the nonwandering set
$\Omega(\widetilde F)=\{x\in\R^k_+$: for every neighborhood $U$ of $x$,
$\widetilde F^n(U) \cap U \neq\varnothing$ for some $n\}$,
$\Omega(\widetilde F)$ is a hyperbolic set and the periodic points are
dense in $\Omega(\widetilde F)$. Results of \citeauthor{przytycki-76}
[(\citeyear{przytycki-76}), 3.11--3.14~and~3.17] imply that key
attributes of Axiom~A endomorphism are: (i) $\Omega(\widetilde F)$
decomposes in a finite number of invariant sets $\Omega^1(\widetilde
F),\ldots, \Omega^m(\widetilde F)$, (ii) for each orbit
$\{x_n\}\subset\Omega^i(\widetilde F)$ of $\widetilde F$, the unstable
manifold at $x_0$ intersects $\Omega^i(\widetilde F)$ in a dense set
and (iii) maps $F$ sufficiently $C^1$ close to $\widetilde F$ are
Axiom~A endomorphisms. Property (iii) implies that $F(x)=D \tilde F(x)$
is an Axiom A endomorphism provided that $D$ is sufficiently close to
the ident ity matrix. Property (ii) implies that each of the invariant
sets $\Omega^i(F)$ is a $\rho$-equivalence class. Linear stability of
the $n^k$ periodic orbits $(g^{n_1}(p),\ldots,g^{n_k}(p))$ with $0\leq
n_j<n$ for $\widetilde F$ implies that, for sufficiently small
perturbations $F$ of $\widetilde F$, $n^k$ of the invariant sets
$\Omega^i(F)$ correspond to linearly stable periodic points, while the
remaining invariant sets are either hyperbolic repellers or
saddles. Since the stable periodic orbits are the only
$\rho$-quasi-attractors,
Theorem~\ref{th.ld} implies the first assertion of the proof. %
%for this part of the proof.

We prove the second assertion. Since $D$ is doubly stochastic, the
nonnegative half-line $\mathcal{L}=\{x\dvtx x_1=\cdots=x_k\geq 0\}$ is
$F$-invariant, that is, $F(x_1 \mathbf{1})=g(x_1) \mathbf{1}$ where
$\mathbf{1}$ is the vector of ones. %and $\circ$ is component wise
%multiplication.
As in the case of the proof of the first assertion, the dynamics of $F$
of restricted $\mathcal{L}$ has the stable periodic point $(p \mathbf
{1}, g(p) \mathbf{1},\ldots, g^{n-1}(p) \mathbf{1})$, %${\mathcal S}
and the complement of its basin of attraction can be decomposed into a
finite number, say $m$, of compact, $g$-invariant sets which have a
dense orbit and are hyperbolic repellers. Define $\widetilde D$ by
$\widetilde d_{ij}=d_{ij}-d_{ik}$ for all $i,j$. By choosing $d_{ij}$
sufficiently close to $1/k$ for all $i,j$, we can make the matrix
$\widetilde D$ as close to zero as we want. Hence,
\citeauthor{earn-levin-06} [(\citeyear{earn-levin-06}), Theorem~1]
implies that $\mathcal{L}$ is a global attractor for the dynamics of
$F$. Moreover, the stable periodic orbit for $F$ restricted to
$\mathcal{L}$ is stable for $F$. Proposition~\ref{prop.dyn} implies
that each of these invariant sets is a $\rho$-equivalence class. Since
the stable periodic orbit is the only $\rho$-quasi-attractor,
Theorem~\ref{th.ld} implies the second assertion.
\end{pf*}

\subsubsection*{Competing species} During the mid twentieth century,
laboratory experiments played a key role in establishing the
competitive exclusion principle in ecology. One classic set of
competition experiments was conducted by \citeauthor{park-48}
(\citeyear{park-48,park-54}) with flour beetles. To model the dynamics of these competing beetles,
collaborators of Park~\citep{leslie-gower-58} used difference
equations, rather than the classical Lotka--Volterra differential
equation model of competition. \citet{cushing-etal-04} showed that
these difference equations exhibit the same dynamical outcomes as the
Lotka--Volterra models. Namely, one or both species may go extinct for
all initial conditions, may coexist about a globally stable equilibrium
or may exhibit contingent exclusion where the initially ``more
abundant'' species excludes the other species. Here, we consider a
stochastic counterpart of the Leslie--Gower model.

Let $N_t^\varepsilon=(N_t^{\varepsilon,1},N_t^{\varepsilon,2})$ and
$X_t^\varepsilon= \varepsilon N_t^\varepsilon$ denote the abundances
and densities of the competing species at time $t$. Once again,
$1/\varepsilon$ corresponds to the volume of their habitat. The
per-capita mean fecundity $f_i(x)$ for species $i$ is given by
$f_i(x)=\frac{b_i}{1+c_{ii}x_i+c_{ij} x_j}$ where $j\neq i$, $b_i>0$ is
the ``intrinsic'' birth rate, $c_{ii}>0$ is the strength of
intraspecific competition and $c_{ij}>0$ is the strength of
interspecific competition. If individual births are independent given
the current density of individuals and Poisson distributed with means
$f_i(X_t^\varepsilon)$ $i=1,2$, then $N_{t}^{\varepsilon,i}$ is a
nonlinear Poisson process associated with the map
$F(x)=(f_1(x)x_1,f_2(x)x_2)$. Proposition~\ref{prop.quasi} implies the
Leslie--Gower process has quasi-stationary distributions
$\mu_\varepsilon$ for $\varepsilon>0$ with $M=\R_+^2=\{x\in\R^2\dvtx
x_i\geq  0\}$ and $M_0=
\partial\R_+^2 =\{x\in\R^2_+\dvtx x_1x_2=0\}$. Results of \citeauthor{cushing-etal-04}
[(\citeyear{cushing-etal-04}), Theorem~4], our Theorems
\ref{th.M0}~and~\ref{th.ld} imply the following result.

%
%th5.4 #&#
\begin{theorem} Let $\mu_\varepsilon$ be a quasi-stationary
distribution for the Leslie--Gower process $\{X_t^\varepsilon\}$. Let
$\mu$ a weak* limit point of these quasi-stationary distributions.
\begin{longlist}
\item[Coexistence.] If $b_i>1$ for $i=1,2$ and $c_i(b_j-1)<b_i-1$
    for $i=1,2$ and $i\neq j$, then $\mu$ is a Dirac measure
    supported by the point
\[
\biggl(\frac{b_2-1}{c_1c_2-1} \biggl(c_1- \frac{b_1-1}{b_2-1} \biggr),
\frac
{b_1-1}{c_1c_2-1} \biggl(c_2- \frac{b_2-1}{b_1-1} \biggr) \biggr).
\]

\item[Extinction or exclusion.] If $b_i<1$ for some $i$, or $b_i>1$
    for $i=1,2$, $b_2-1>(b_1-1)/c_1$ and $b_1-1>(b_2-1)/c_2$, or
    $b_i>1$ for $i=1,2$, $b_2-1<(b_1-1)/c_1$ and
    $b_1-1<(b_2-1)/c_2$, then $\mu$ is supported by $\partial
    \R_+^2$.
\end{longlist}
\end{theorem}

The case for which our results are not conclusive is when the dynamics
of the Leslie--Gower model are bistable [i.e., $b_i>1$ for $i=1,2$ and
$c_i(b_j-1)>b_i-1$ for $i=1,2$ and $i\neq j$] in which case there is a
positive unstable equilibrium and all initial conditions not lying on
its stable manifold (which has dimension one) go to~$\partial\R_+^2$.
However, we conjecture that $\mu$ is supported on the boundary of the
positive quadrant in this case.

\subsubsection*{Host-parasitoid interactions} Predator-prey
interactions involve one\break  species benefiting by harming another species.
These interactions are the fundamental\vadjust{\goodbreak} building blocks for all food
webs. An important class of predators is parasitoids such as wasps or
flies whose young develop in and ultimately kill their host~\citep
{godfray-94}. Mathematical models of these interactions have been
studied for almost a century [\citet{thompson-24},
\citet{nicholson-bailey-35}, \citeauthor{hassell-78} (\citeyear{hassell-78,hassell-00}),
\citet{may-95}, \citeauthor{jmb-06} (\citeyear{jmb-06,jbd-07}),
\citet{gidea-etal-11}].
As predator-prey interactions are inherently oscillatory, these studies
often focused on identifying mechanisms that stabilize predator-prey
interactions.

Here, we introduce a stochastic analog of these deterministic models.
Let $H_t^\varepsilon$ and $P_t^\varepsilon$ denote the abundances of
host and the parasitoid in generation $t$, respectively. Let
$X_t^{\varepsilon}=\varepsilon N_t^{\varepsilon} = \varepsilon
(H_t^\varepsilon, P_t^\varepsilon)$ be their densities where
$1/\varepsilon$ is the size of the environment. Let
$f(X_t^{\varepsilon,1})$ be the mean number of offspring produced by
an individual host. Let $g(X_t^{\varepsilon})$ be the probability that
an offspring escapes parasitism from the parasitoids. We update the
population state $X_t^\varepsilon$ according to the following rules:
\begin{itemize}
\item Each adult host independently produces a Poisson distributed
number of offspring with mean $f(X_t^{\varepsilon,1})$. Let
$M_{t+1}$ be the total number of offspring which is Poisson
distributed with mean $H_{t}^\varepsilon
f(X_t^{\varepsilon,1})$.

\item Each offspring survives parasitism independently with
probability $g(X_t^\varepsilon)$. Let $H_{t+1}^\varepsilon$
equal the number of surviving offspring and
$P_{t+1}^\varepsilon=M_{t+1}-H_{t+1}^\varepsilon$ be the number
of parasitized offspring which all emerge as parasitoids in the
next generation.
\end{itemize}
Since $(H_{t+1}^\varepsilon,P_{t+1}^\varepsilon)$ is binomial
distributed with $M_{t+1}$ trials,
$(H_{t+1}^\varepsilon,P_{t+1}^\varepsilon)$ are independent Poisson
random variables. Hence, $X^\varepsilon_t$ is a nonlinear Poisson
process associated with the map
\[
F(x)=\bigl(f(x_1)x_1g(x), f(x_1)x_1
\bigl(1-g(x)\bigr)\bigr)
\]
on $\R^2_+=\{(x_1,x_2)\in\R^2 \dvtx x_i \geq 0\}$. Provided that $F$ is
continuous, and $f$ is a compact map, Proposition~\ref{prop.quasi}
implies that there is a quasi-stationary distribution $\mu_\varepsilon$
for $X_t^\varepsilon$ with $\varepsilon>0$.

To understand the support of the weak* limit points $\mu$ of
$\mu_\varepsilon$, we focus on a generalized Thompson
model [\citet{thompson-24}, \citet{getz-mills-96},
\citeauthor{jmb-06} (\citeyear{jmb-06,jbd-07})]. For this model,
$f(x_1)=\exp(r(1-x_1/K))$ is given by the Ricker equation where $r>0$
is the intrinsic rate of growth of the host, and \mbox{$K>0$} is the
host's carrying capacity. The escape function $g(x)=(1+x_2/(bx_1
k))^{-k}$ corresponds to a negative binomial escape function with
egg-limited encounter rates. Here, $b>0$ is the attack rate of the
parasitoid and $1/k>0$ represents how ``clumped'' or ``aggregated''
parasitoid attack are; that is, smaller $k$ correspond to greater
aggregation of parasitoid attacks. Notice that while $g$ is not defined
at $x_1=0$, the map $F$ extends continuously to $x_1=0$ if we set
$F(x)=0$ whenever $x_1=0$. Combining results from \citeauthor{jbd-07}
[(\citeyear{jbd-07}), Theorem~3.1, 3.2] and Theorem~\ref{th.ld} yields
the following results for $k<1$, that is, parasitoid attacks are
sufficiently aggregated~(\citet{hassell91}).

%
%th5.5 #&#
\begin{theorem} Let $\mu$ be a weak* limit point of the
quasi-stationary distributions for the Thompson host-parasitoid process
$X_t^\varepsilon$. Assume $k<1$, and define
\[
y^*=\max\bigl\{y\geq 0\dvtx \exp(-r) \bigl(\bigl(1+y/(bk)\bigr)^{k}-1
\bigr)=y\bigr\}.
\]
Then:
\begin{longlist}
\item[Extinction.] If $\exp(r)(1+y^*/(bk))^{-k}<1$, then $\mu$ is
    supported by the $\partial\R^2_+=\{x\in\R^2_+\dvtx x_1x_2=0\}$.
\item[Coexistence.] If $\exp(r)(1+y^*/(bk))^{-k}>1$, then $\mu$ is
    supported by \mbox{$\R^2_+\setminus\partial\R^2_+$}. Moreover, for an
    open and dense set of parameter values $(r,b)$ satisfying
    $\exp(r)(1+y^*/(bk))^{-k}>1$, $\mu$ is supported by a periodic
    orbit.
\end{longlist}
\end{theorem}
When $k\geq 1$ (i.e., parasitoid attacks are not sufficiently
aggregated),\break \citeauthor{jbd-07}  [(\citeyear{jbd-07}), Theorem~3.1]
implies coexistence does not occur for the deterministic model.
However, this extinction often involves unstable sets in the interior
of $\R^2_+$. Consequently, our results are not applicable. Nonetheless,
we conjecture that $\mu$ is supported by $\partial\R^2_+$.

%%%%%%%%%%%%%%%%%%%%%%%%%%%%%%%%%%%%%%%%%%%%%%%%%%%%%%%%%%%%%%%%%%%%%%%%%%%%%%%%%%%%%%%%%%%%%%%
%s5.3 #&#
%s5.3 ###
\subsection{Multinomial processes}\label{s5.3} Consider a landscape
with $N$ sites that can be in one of $k$ states. These states may
correspond to occupation by individuals playing different strategies in
the context of evolutionary games, or different genotypes in the
context of population genetics. Let $M=\{x\in\R^k\dvtx x_i \geq 0,
\sum_{i=1}^k x_i=1\}$ be the $k$-simplex and $F\dvtx M\to M$ be a
continuous map.

For each $x\in M$, let $Z_1(x),Z_2(x),\ldots$ be a sequence of
independent random vectors with a multinomial distribution with $N$
trials, $k$ possible outcomes and probability $x_i$ of producing type
$i$ in a single trial. If $\varepsilon=1/N$ and $X_0^\varepsilon\in
M\cap\varepsilon\mathbb{Z}^k$ is given, then we can define a Markov
chain $\{X_t^\varepsilon\}_{t=0}^\infty$ on $M\cap\varepsilon
\mathbb{Z}^k$ iteratively by
\[
X^\varepsilon_{t+1}=\varepsilon Z_{t+1}\bigl(F
\bigl(X_t^{\varepsilon}\bigr)\bigr).
\]
Since $\{X_t^\varepsilon\}$ is a finite-state Markov chain,
quasi-stationary distributions exist uniquely whenever the transition
matrix restricted to the transient states is aperiodic and
irreducible~\citep{DarSen65}. When $M_0$ is the boundary of the simplex
and $F_i(x) = x_i f_i(x)$, with $f_i$ continuous and positive, we
therefore always have a unique quasi-stationary distribution, and we
prove in Proposition~\ref{pr.largedev2} of Section~\ref{s6} that the
large deviation hypotheses are satisfied. As a particular application
of these multinomial processes, we consider evolutionary games.

\subsubsection*{Evolutionary game dynamics}
Evolutionary game theory studies the dynamics of populations of
players, each programmed to play a fixed strategy throughout their life
time [\citeauthor{hofbauer-sigmund-98}
(\citeyear{hofbauer-sigmund-98,hofbauer-sigmund-03}),
\citet{cressman-03}]. These populations often exhibit
frequency dependent selection; the reproductive success of a player
changes in time due to the composition of strategies in the population.
The study of evolutionary games has led to fundamental insights into
the evolution of animal conflicts~\citep{maynardsmith-74},
cooperation~\citep{nowak-etal-04,imhof-nowak-10}, habitat
selection~\citep{schreiber-fox-getz-00,cressman-etal-04} and mating
systems~\citep{sinervo-lively-96}.

A basic deterministic model for evolutionary games is the discrete-time
replicator equation~\citep{hofbauer-sigmund-03},
\[
F_i(x)= x_i \frac{\sum_{j}a_{ij}x_j+c}{\sum_{jl} a_{jl}x_jx_l+c},
\]
where $x_i$ is the frequency of strategy $i$ in the population, the
entries $a_{ij}$ of the ``pay-off'' matrix $A$ describe the fitness
gain to strategy $i$ when interacting with strategy~$j$, and $c$ is the
``basal'' fitness of an individual. The dynamics of this discrete
system have been studied
extensively [\citeauthor{hofbauer-sigmund-98} (\citeyear{hofbauer-sigmund-98,hofbauer-sigmund-03})]. Here, we
describe a stochastic analog of these games that account for finite
population sizes. In the case of two-strategies, this stochastic analog
corresponds to the frequency dependent Wright--Fisher processes studied
by \citet{imhof-nowak-06}.

Let $N_t^{\varepsilon, i}$ and $X_t^{\varepsilon,i}$ denote the
abundance and frequency of the $i$th strategy at time~$t$. Here,
$1/\varepsilon$ is assumed to be an integer that corresponds to the
total population size which does not change over time. If we update
$N_t^\varepsilon$ by taking a multinomial random variable with
$1/\varepsilon$ trials and probabilities $F(X_t^\varepsilon)$, then we
get a multinomial process. One can interpret this process as
individuals producing many offspring proportional to their fitness
$\sum_j a_{ij} x_j +c$, and randomly selecting individuals from the
offspring ``pool'' to replace their parents. For this stochastic
process $\{X_t^\varepsilon\}$, $M_0 = \{ x\in M\dvtx \prod x_i =0\}$ is
an absorbing set that corresponds to the loss of one or more
strategies.

We can leverage two results from the theory of replicator dynamics to
describe the support of the quasi-stationary distributions
$\mu_\varepsilon$ when $\varepsilon>0$ is sufficiently small, and the
basal payoff $c$ is sufficiently large. To first order in $1/c$,
\[
F_i(x)\approx x_i + x_i \frac{1}{c}
\biggl(\sum_{j}a_{ij}x_j-
\sum_{jl} a_{jl}x_jx_l
\biggr)
\]
and the dynamics of the map $x\mapsto F(x)$ can be viewed as a
Cauchy--Euler approximation to the classical continuous time replicator
equations
%
%e11 #&#
%e10 ###
\begin{equation}
\label{eq.replicator} \frac{dx_i}{dt} = x_i \biggl(\sum
_{j}a_{ij}x_j-\sum
_{jl} a_{jl}x_jx_l
\biggr).
\end{equation}
Using this observation and work of
\citet{hofbauer-sigmund-98,garay-hofbauer-03}, we can prove a
sufficient condition for the stochastic replicator processes exhibiting
metastable persistence for $\varepsilon>0$ small and $c>0$ large.

%
%th5.6 #&#
\begin{theorem}\label{thm.replicator} Let $\mu$ be a weak* limit point
of the quasi-stationary distributions $\mu^\varepsilon$ for the
replicator process with $c>0$ sufficiently large.
If there exists $p_i>0$ such that
\[
\sum_i p_i \biggl(\sum
_{j}a_{ij}x_j-\sum
_{jl} a_{jl}x_jx_l
\biggr)>0
\]
at every equilibrium $x\in M_0$ for $F$, then $\mu$ is supported by
$M_1$.
\end{theorem}

\begin{pf}
\citeauthor{hofbauer-sigmund-98} [(\citeyear{hofbauer-sigmund-98}),
Theorem~13.6.1] implies that the continuous-time replicator equations
\eqref{eq.replicator} have a global attractor $A\subset M_1$ whenever
there positive weights $p_1,\ldots, p_k$ such that
\[
\sum_i p_i \biggl(\sum
_{j}a_{ij}x_j-\sum
_{jl} a_{jl}x_jx_l
\biggr)>0
\]
at every equilibrium $x\in M_0$ for $F$. In particular, in the
terminology of \citeauthor{garay-hofbauer-03}
[(\citeyear{garay-hofbauer-03}), Definition~2.1], \eqref{eq.replicator}
admits a good average Lyapunov function whenever these positive weights
exist. Consequently, \citeauthor{garay-hofbauer-03}
[(\citeyear{garay-hofbauer-03}), Theorem~8.3] implies that the
discrete-time replicator equation $F$ admits a global attractor
$A\subset M_1$ whenever $c$ is sufficiently large. This global
attractor, however, need not be a unique $\rho$-quasi attractor.
However, as suggested by Comment 3.8, we can apply
Lemma~\ref{lm.lambda}(b$'$) to conclude that $\mu(V_0)=0$ for some
neighborhood $V_0$ of $M_0$.
\end{pf}

As an interesting special case consider the rock-paper-scissor game
where the payoff-matrix is of the form
\[
A=\pmatrix{ 0&-a_2 & b_3
\cr
b_1&0&
-a_3
\cr
-a_1&b_2&0}
\]
with $a_i$ and $b_i$ positive. \citet{zeeman-80} proved that if
$\mathrm{det}(A)>0$, then the persistence condition of
Theorem~\ref{thm.replicator} is satisfied. Moreover, for the
continuous-time replicator equations, there is a globally stable
internal equilibrium. For $c>0$ sufficiently large, this equilibrium is
also globally stable for the map $F$ and, consequently $\mu$ is a Dirac
measure supported by this equilibrium. When $\mathrm{det}(A)<0$, Zeeman
proved that the internal equilibrium is unstable and all other orbits
of the continuous-time deterministic system approach the boundary and
one can show that the same conclusion holds for the discrete-time
system when $c>0$ is sufficiently large. Since the boundary is not a
global attractor in this case, we cannot apply Theorem~\ref{th.M0}.
None the less, we conjecture that $\mu$ is supported on the boundary
$M_0$ in this case.

%%%%%%%%%%%%%%%%%%%%%%%%%%%%%%%%%%%%%%%%%%%%%%%%%%%%%%%%%%%%%%%%%%%%%%%%%%%%%%%%%%
%%%%%%%%%%%%%%%%%%%%%%%%%%%%%%%%%%%%%%%%%%%%%%%%%%%%%%%%%%%%%%%%%%%%%%%%%%%%%%%%%%
%s6 #&#
%s6 ###
\section{Large deviation results for Poisson and multinomial
models}\label{s6}

In this section we prove the existence of quasi-stationary
distributions (as needed) and verify our large deviation hypotheses for
the Poisson and multinomial models introduced in Section~\ref{s5}.

%s6.1 #&#
%s6.1 ###
\subsection{Nonlinear Poisson branching model}\label{s6.1}

We first prove the existence of quasi-stationary distributions for the
nonlinear Poisson processes introduced in Section~\ref{s5.2}.

%
%pr6.1 #&#
\begin{proposition}\label{prop.quasi} If $\sup_{x\in\R^k_+} \| F(x)\|
< \infty$ and $F_i(x)>0$ for all $x\in M\setminus M_0$ and $i$, then
the nonlinear Poisson process $\{X_t^\varepsilon\}$ associated with $F$
has at least one quasi-stationary distribution supported on $M\setminus M_0$.
\end{proposition}

\begin{pf} For notational convenience, we prove this result for the Markov
chain $\{N_t^\varepsilon\}$ and call $p^{\varepsilon}(x,y)$ its
transition kernel. Let $X=\mathbb{Z}^k_+ \setminus M_0$ and
$q^\varepsilon(x,y)$ denote the restriction of $p^{\varepsilon}$ to
$X$. Let $l^1(X)$ denote absolutely summable functions from $X$ to
$\R$. For $u\in l^1(X)$, define $\|u\|_1=\sum_x |u(x)|$.

We can define a linear operator $Q^\varepsilon$ from $l^1(X)\to l^1(X)$
by $(uQ^\varepsilon)(x)=\sum_{y \in X} q^\varepsilon(y,x)u(y)$. Recall
that $\mu_{\varepsilon}$ is a QSD for the Markov chain
$N^{\varepsilon}$ if and only if it is a nonnegative eigenvector of the
operator $Q^{\varepsilon}$. Since $F_i(x)>0$ for all $i$ and $x\in
M\setminus M_0$, $q^\varepsilon(x,y)>0$ for all $x,y\in X$.
%Let $Q^\varepsilon\in\mathcal{L}(l^1(X))$ be the operator induced by
%the nonnegative kernel $q^{\varepsilon}$.

Next, we show $Q^\varepsilon$ is a compact operator, that is, the image
of the unit ball under $Q^\varepsilon$ is precompact. Recall that a
closed set $C$ in $l^1(X)$ is compact if and only if it is bounded and
equisummable: for all $\delta>0$, there exists $N$ such that $\sup_{u
\in C} \sum_{\|x\| \geq N} |u(x)| \leq\delta$. Defining
$G(x)=F(\varepsilon x)/\varepsilon$, we get
\[
q^\varepsilon(x,y)=\prod_{i=1}^k
\frac{G_i(x)^{y_i}}{y_i!} \exp \bigl(-G_i(x)\bigr)\leq \prod
_{i=1}^k \frac{ m^{y_i} }{ y_i!}\qquad\mbox{where } m=\sup
_{x\in X} \bigl\|G(x)\bigr\|.
\]
Hence, for $u$ with $\|u\|_1\leq 1$,
\begin{eqnarray*}
\bigl\|uQ^\varepsilon\bigr\|_1&=&\sum_x\Biggl|
\sum_y q^\varepsilon(y,x)u(y)\Biggr|
\\
&\leq & \sum_x \sum_y
\prod_{i=1}^k \frac{m^{x_i} }{ x_i!} \bigl|u(y)\bigr|
\\
&\leq & \sum_x \prod_{i=1}^k
\frac{m^{x_i} }{ x_i!}=e^{km}.
\end{eqnarray*}
Moreover, given $\delta>0$, choose $N>0$ such that $\sum_{\|x\| \geq N}
\prod_{i=1}^k m^{x_i}/x_i! < \delta$. Then $\sum_{|x|\geq  N}
|(uQ^\varepsilon)(x)| \leq \delta$ for all $u$ such that $\|u \|_1\leq
1$. Hence, $Q^\varepsilon$ is a compact operator.

On the other hand, given that $q^{\varepsilon}$ is strictly positive on
$X \times X$, we have for any $Y \subset X$,
\[
\sum_{x \notin Y} \sum_{y \in Y}
q^{\varepsilon}(x,y) >0.
\]
Applying the following result of Jentzsch on Kernel positive operators
completes the proof of this proposition.
\end{pf}
%
%
%th6.2 #&#
\begin{theorem}
[{[\citeauthor{Schaefer74} (\citeyear{Schaefer74}),
Theorem~\textup{V}.6.6]}] Let $E = L^p(\mu)$, where $1 \leq p \leq+
\infty$ and $(X,\Sigma,\mu)$ is a $\sigma$-finite measure space.
Suppose $Q \in\mathcal{L}(E)$ is an operator given by a
$(\Sigma\times\Sigma)$-measurable kernel $q \geq0$, satisfying the
following two assumptions:
\begin{longlist}[(ii)]
\item[(i)] some power of $Q$ is compact;

\item[(ii)] $Y \subset\Sigma$, $\mu(Y) >0$ and $\mu(\Sigma \setminus
    Y)>0$ implies
\[
\int_{X \setminus Y} \int_{Y} q(x,y) \mu(dx)
\mu(dy) >0.
\]
\end{longlist}
Then the spectral radius $r(Q)$ is positive, is a simple eigenvalue,
its unique renormalized eigenvector $v$ satisfies $v(x)>0$, $\mu$
almost surely and $r(Q)$ is the only eigenvalue of $Q$ with a positive
eigenvector. Moreover, if $q(x,y)>0$ $(\mu\otimes\mu)$ almost surely,
then every other eigenvalue of $Q$ has modulus strictly smaller
than~$r(Q)$.
\end{theorem}

The following proposition verifies the large deviation hypotheses of
Section~\ref{s2} for the nonlinear Poisson processes.

%
%pr6.3 #&#
\begin{proposition} \label{pr.largedev} Assume that $x \mapsto F(x)$ is
continuous, $F_i(x)>0$ for all $i$ and $x\in M\setminus M_0$, and
$\sup_{x\in\R^k_+} \| F(x)\| < \infty$. Then the nonlinear Poisson
process $\{X_t^\varepsilon\}$ associated with $F$ satisfies
Hypotheses~\ref{hy.rho} and \ref{hy.rho2}.
\end{proposition}

\begin{pf} Let $\mu^x_{\varepsilon}$ denote the distribution of the Poisson
random vector $X_{t+1}^{\varepsilon}$ conditional to
$X^{\varepsilon}_t = x$. The logarithmic moment generating function
relative to $\mu^x_{\varepsilon}$ is given by
\[
\Lambda_{\varepsilon,x}(\lambda) = \log\mathbb{E} \bigl(e^{\langle \lambda,
\varepsilon Z_1(F(x)/\varepsilon)\rangle} \bigr)
= \sum_{i=1}^k \frac
{F_i(x)}{\varepsilon}
\bigl(e^{\varepsilon\lambda_i} - 1\bigr).
\]
Hence, the family $\varepsilon
\Lambda_{\varepsilon,x}(\cdot/\varepsilon)$ is identically equal on
$\R^k$ to the function $\Lambda_{x}(\lambda) = \sum_{i=1}^k F_i(x)
(e^{\lambda_i} -1)$. Thus, by the G\"artner--Ellis theorem [see, e.g.,
\citet{dembo-zeitouni-93}, Theorem~2.3.6], the family
$\mu^x_{\varepsilon}$ satisfies a large deviation principle with convex
rate function $\Lambda_x^*(y)=\sum_{i=1}^k y_i \log\frac{y_i}{F_i(x)} +
F_i(x) - y_i$; that is, for any closed set $F\subset\R^k_+$ and
$x\in\R^k_+$,
\[
\limsup_{\varepsilon\to0} \varepsilon\log\mu_{\varepsilon}^x(F)
\leq  - \inf_{y\in F} \Lambda_x^*(y)
\]
and for any open set $G\subset\R^k_+$,
\[
\liminf_{\varepsilon\to0 } \varepsilon\log\mu_{\varepsilon}^x(G)
\geq - \inf_{y\in G} \Lambda_x^*(y).
\]
Hence, if we define
\[
\rho(x,y)=\Lambda^*_x(y),
\]
then $\rho$ immediately satisfies (i), (ii) and the upper bound of (iv)
of Hypothesis~\ref{hy.rho}.

Now let us derive the uniform lower bound of Hypothesis~\ref{hy.rho}.
Pick a compact set $K \subset M_1$ and an open ball $B \subset M$. Let
$B^{\varepsilon}=B\cap(\varepsilon\mathbb{N}^k)$ be the
$\varepsilon$-lattice on $B$. Let $\alpha>0$. For every $x \in K$,
there exists $y(x) \in B$ such that $\rho(x,y(x)) \leq\inf_{y \in B}
\rho(x,y) + \alpha$. Choose $\varepsilon_0>0$ small enough such that
\[
d\bigl(y,y'\bigr)< \varepsilon_0\quad\Rightarrow\quad\bigl|
\rho(x,y)-\rho\bigl(x,y'\bigr)\bigr| < \alpha
\]
for all $y,y'\in B$ and $x\in K$. For each $x \in K$ and each
$0<\varepsilon< \varepsilon_0$, we can choose a point
$y_{\varepsilon}(x) = \varepsilon
(n_1^{\varepsilon},\ldots,n_k^{\varepsilon})$ such that
$d(y_\varepsilon(x),y(x))<\varepsilon$ and, consequently,
$\rho(x,y_{\varepsilon}(x)) \leq\inf_{y \in B} \rho(x,y) + 2 \alpha$.
For all $0<\varepsilon< \varepsilon_0$ and all $x \in K$, we have
\[
\mu_{\varepsilon}^x (B) \geq\mu_{\varepsilon}^x
\bigl(\bigl\{y_{\varepsilon
}(x)\bigr\} \bigr) = \prod_{i=1}^k
e^{-F_i(x)/\varepsilon} \frac
{(F_i(x)/\varepsilon)^{n_i^{\varepsilon}}}{n_i^{\varepsilon}!}.
\]
Recall that, for any $p \in\mathbb{N}$, $-p-\log p! + p \log p \geq
-(1+\log p)$ and define $I_+=\{i\in \{1,\ldots,k\}\dvtx
n_i^\varepsilon>0\}$. A straightforward computation gives
\begin{eqnarray*}
\varepsilon\log\mu_{\varepsilon}^x (B) &\geq&- \sum
_{i=1}^k F_i(x) +\sum
_{i\in I_+} \biggl( \varepsilon n_i^{\varepsilon} \log
\frac
{F_i(x)}{\varepsilon} - \varepsilon\log\bigl(n_i^{\varepsilon}!\bigr)
\biggr)
\\
&=& - \rho\bigl(x,y_{\varepsilon}(x)\bigr) + \sum_{i\in I_+}
\varepsilon n_i^{\varepsilon} \biggl(-1 - \frac{1}{n_i^{\varepsilon}} \log
\bigl(n_i^{\varepsilon}!\bigr) + \log n_i^{\varepsilon}
\biggr)
\\
&\geq& -2 \alpha- \inf_{y \in B} \rho(x,y) - \varepsilon\sum
_{i\in
I_+} \bigl(1 + \log n_i^{\varepsilon}
\bigr).
\end{eqnarray*}
The last quantity goes to zero as $\varepsilon$ goes to zero,
independently of $x$ since $n_i^{\varepsilon}$ is of order
$(y_{\varepsilon}(x))_i/\varepsilon$ and the quantities
$(y_{\varepsilon}(x))_i$ are bounded. Hence, we have shown that the
lower bound (\ref{lower}) holds uniformly for any compact set $K\subset
M_1$ and open ball $B$.

To verify that $\rho$ satisfies (iii) of Hypothesis~\ref{hy.rho}, we
first prove the following lemma.

%
%le6.4 #&#
\begin{lemma} \label{lm.g} Define $g\dvtx (0,\infty)\times[0,\infty)\to
[0,\infty)$ by $g(x,y)=y\log\frac{y}{x}+ x-y$. Then for all
$\delta>0$ and $m>0$,
\[
\inf\bigl\{g(x,y)\dvtx |x-y|\geq \delta, x\leq  m\bigr\}\geq a >0.
\]
\end{lemma}

\begin{pf} Let $m>0$ and $\delta>0$. We have that $C:=\{(x,y)\dvtx|x-y|\geq
\delta, x\leq  m\}=A\cup B$ with $A=\{(x,y)\dvtx m\geq  x \geq \delta,
0\leq  y \leq  x-\delta\}$ and $B=\{(x,y)\dvtx y \geq \delta,
0<x\leq y-\delta\}$. Since $A$ is compact and $g$ restricted to $A$ is
positive and continuous, $\inf_{(x,y)\in A} g(x,y)>0$. Restricted to
$B$, $g(x,y)$ is positive and increasing in $y$. Hence, $\inf_{(x,y)\in
B} g(x,y) = \inf_{0<x\leq  m} g(x,x+\delta)$. Since
\[
\frac{d}{dx} g(x,x+\delta) = \log(1+\delta/x)-\delta/x <0,
\]
we get $\inf_{0<x\leq  m} g(x,x+\delta)= g(m,m+\delta)>0$. Thus,
$\inf_{(x,y)\in C} g(x,y)>0$.
\end{pf}

Let now $\delta>0$ be given, $d(x,y)=\max_i |x_i-y_i|$ and
$g$ be as defined in the lemma. If $|y_i-F_i(x)|\geq \delta$, then
Lemma \ref{lm.g}, with $m=\sup_{x\in\R^k_+} \|F(x)\|$, implies that
$\rho(x,y)=\sum_j g(F_j(x),y_j)\geq  g(F_i(x),y_i)\geq  a$.

To check that the uniform upper bound \eqref{upper} holds, notice that
it is sufficient to prove that quantities $
\mu_{\varepsilon,i}^x([F_i(x)+ \delta, + \infty[)$ (where
$\mu_{\varepsilon,i}^{x}(\cdot)$ is the $i$-marginal
of~$\mu^{x}_{\varepsilon}$, namely the distribution of the $i$
component of $X_{t+1}^{\varepsilon}$, conditional to
$X_t^{\varepsilon}=x$) are bounded above by some expression which goes
to zero as $\varepsilon$ goes to zero, uniformly in $x \in K$. This is
an easy consequence of Chernov's upper bound,
\[
\mu_{\varepsilon,i}^x\bigl([F_i(x)+ \delta, + \infty[
\bigr) \leq e^{-(1/\varepsilon) g(F_i(x),F_i(x)+\delta)} \leq e^{-1/\varepsilon \beta},
\]
where $\beta= \inf \{g(x,y) \dvtx |x-y|\geq\delta, x \leq m \} >0$, and
the quantity on the right-hand side goes to zero uniformly in $x$
by~(iii).

Finally to verify Hypothesis~\ref{hy.rho2}, notice that
$p^\varepsilon(x,0)=\exp(-\sum_i F_i (x)/\varepsilon)$, and recall that
$F(0)=0$. Hence, given $c>0$, choose a neighborhood $V_0$ of $\{0\}$
such that $\sum_i F_i(x)\leq  c$ whenever $x\in V_0$. Then $\varepsilon
\log p^\varepsilon(x,0) \geq -c $ whenever $x\in V_0$.
\end{pf}

%%%%%%%%%%%%%%%%%%%%%%%%%%%%%%%%%%%%%%%%%%%%%%%%%%%%%%%%%%%%%%%%%%%%%%%%%%%%%%%%%%%%%%%%%%%%%%%%%
%s6.2 #&#
%s6.2 ###
\subsection{Multinomial model}\label{s6.2}

Here we verify the large deviation assumptions for the multinomial
processes introduced in Section~\ref{s5.3}.

%
%pr6.5 #&#
\begin{proposition} \label{pr.largedev2} Assume $F_i(x)=x_i f_i(x)$,
with $f_i$ continuous and positive. Then the multinomial process $\{
X_t^\varepsilon\}$ associated with $F$ satisfies
Hypotheses~\ref{hy.rho}~and~\ref{hy.rho2} with respect to the absorbing
set $M_0=\{ x\in M\dvtx \prod x_i=0\}$.
\end{proposition}

\begin{pf} Let $\mu_N^x$ be the law of the multinomial random vector
$\frac{1}{N} Z_1(F(x))$, which can be written as $\frac{1}{N}
\sum_{i=1}^N Y_i(F(x))$, where $(Y_i(F(x)))_i$ is an i.i.d. sequence with
distribution $\mathbb{P}[Y_1(F(x)) = e_j] = F_j(x)$ ($e_j$ is the
unitary vector in $\mathbb{R}^k$ whose $j$th component equals one).
By\vadjust{\goodbreak}
Cram\'er's theorem [see, e.g., \citet{dembo-zeitouni-93},
Theorem~2.2.30], the sequence $\mu_N^x$ satisfies a large deviation
principle with convex rate function $\Lambda_x^*(y)=\sum_{i=1}^k y_i
\log\frac{y_i}{F_i(x)}$. Hence, if we define $\rho(x,y) =
\Lambda^*_x(y)$, then $\rho$ immediately satisfies (i), (ii), (iv) of
Hypothesis~\ref{hy.rho}.

The proof of the uniform lower bound is similar to the Poisson
branching process case. Let $K\subset M_1$ be a compact set and
$B\subset M$ an open ball. Let $B^N = B\cap\frac{1}{N} \mathbb{N}^k$.
Let $\alpha>0$ be given. For every $x \in K$, there exists $y(x)\in B$
such that $\rho(x,y(x))\leq \inf_{y\in B} \rho(x,y)+\alpha$. Choose
$N_0\geq 1$ sufficiently large such that
\[
d\bigl(y,y'\bigr)<1/N_0\quad\Rightarrow\quad\bigl|\rho(x,y)-\rho
\bigl(x,y'\bigr)\bigr|<\alpha
\]
for all $x\in K$ and $y,y'\in B$. For each $x\in K$ and $ N\geq  N_0$,
we choose $y_{N}(x) = \frac{1}{N} (n_1^N,\ldots,n_k^N)$ such that
$d(y_N(x),y(x))<1/N$. Let $I_+=\{i\dvtx n^N_i >0\}$. For $N$ large
enough,
\begin{eqnarray*}
\frac{1}{N} \log\mu_N^x(B) &\geq&
\frac{1}{N} \log\mu_N^x\bigl(\bigl
\{y_N(x)\bigr\}\bigr)
\\[-2pt]
&=& \frac{1}{N} \log \biggl(N!\prod_{i=1}^k
\frac{(F_i(x))^{n_i^N}}{
n_i^N!} \biggr)
\\[-2pt]
&\geq& - \rho\bigl(x,y_N(x)\bigr) + \biggl(1 + \frac{1}{N}
\log N! - \log N \biggr)
\\[-2pt]
&&{} + \frac{1}{N} \sum_{i\in I_+}
n_i^N \biggl(-1 - \frac{1}{n_i^N}\log
n_i^N! + \log n_i^N \biggr)
\\[-2pt]
&\geq& -2 \alpha- \inf_{y \in B} \rho(x,y) - \frac{1}{N}
\sum_{i\in
I_+} \bigl(1 + \log n_i^N
\bigr).
\end{eqnarray*}
The last quantity goes to zero as $N\to\infty$, independently of $x$
since $n_i^{N}$ is of order $N(y_{N}(x))_i$ and the quantities
$(y_{N}(x))_i$ are bounded. Hence, we have shown that the lower bound
(\ref{lower}) holds uniformly for any compact set $K\subset M_1$ and
open ball~$B$.

To verify that $\rho$ satisfies (iii) of Hypothesis~\ref{hy.rho}, assume by
contradiction that there exist $\beta>0$ and two sequences $(x_n)_n$
and $(y_n)_n$ in the $k$-simplex $M$, converging, respectively, to $x$
and $y$, and such that
\[
\lim_n \rho(x_n,y_n) = 0 \quad\mbox{and}\quad d(x_n,y_n) \geq\beta.
\]
Define\vspace*{-4pt} $I_0 = \{i \in\{1,\ldots,k\} \dvtx x^i = 0\}$.
Notice that, if $i \in I_0$, then $\lim_n y_n^i \log\frac{y_n^i}{x^i_n}
= 0$ if $y^i = 0$ and $y_n^i \log\frac{y_n^i}{x^i_n} \rightarrow+
\infty$ otherwise. As a consequence, $y^i = 0\ \forall i \in I_0$. We
consider two cases separately.

Assume first that all the components of $x$ are zero except the first
one, $I_0 = \{2,\ldots,k\}$. Then $y^i=0$ for $i=2,\ldots,k$, which
implies that $x=y$, a contradiction.\vadjust{\goodbreak}

Assume now that $I_0$ contains at most $k-2$ terms. Call $I_1$ its
complementary and assume without loss of generality that $I_1 =
\{1,\ldots,n_1\}$. Define $\tilde{x} = \{x^1,\ldots,x^{n_1}\}$,
$\tilde{y} = \{y^1,\ldots,y^{n_1}\}$ and notice that $\tilde{x} \in
\operatorname{Int}(\Delta^{n_1})$ and $\tilde{y} \in \Delta^{n_1}$,
where $\Delta^{n_1}$ is the \mbox{$n_1$-}simplex. Define analogously
the sequences $(\tilde{x}_n)_n$ and $(\tilde{y}_n)_n$, which belong to
the set $\{u \in\R^{n_1}\dvtx u^i \geq0 \sum_i u^i \leq1\}$. Let now
$\tilde{\rho}$ be the application given by
\[
(u,v) \in\mathbb{R}^{n_1}_* \times\mathbb{R}^{n_1} \mapsto\sum
_{i=1}^{n_1} v_i \log
\frac{v_1}{u_i}.
\]
This map is continuous and strictly positive in
$(\tilde{x},\tilde{y})$, since $|\tilde{x} - \tilde{y}|>\beta$.
Therefore, there exists $\delta>0$ and $r>0$ such that $\tilde{\rho} >
\delta$ on $B_{\R^{n_1}}(\tilde{x},r) \times
B_{\R^{n_1}}(\tilde{y},r)$. This\vadjust{\goodbreak} concludes the
proof since
\begin{eqnarray*}
\liminf_{n \rightarrow+ \infty}
\rho(x_n,y_n) &=& \liminf_{n
\rightarrow+ \infty} \biggl(
\tilde{\rho}(\tilde{x}_n,\tilde{y_n}) + \sum
_{i \in I_0} y^i_n \log\frac{y^i_n}{x^i_n}
\biggr)
\\[-2pt]
&\geq& \beta+ \lim_{n \rightarrow+ \infty} \sum_{i \in I_0}
y^i_n \log\frac{y^i_n}{x^i_n} = \beta.
\end{eqnarray*}

The uniform upper bound \eqref{upper} holds by an application of
Chernov upper bound. Let $\delta>0$. Then
\[
\mu_{N,i}^x\bigl(\bigl[F_i(x)+\delta,+ \infty\bigr[
\bigr) \leq e^{-N\beta},
\]
where $\beta= \inf \{y \log\frac{y}{x}, (x,y) \in(0,1]^2,
|x-y| \geq\delta \} >0$.

Finally to verify Hypothesis~\ref{hy.rho2}, we have
$p^\varepsilon(x,M_0)\geq \max_i (1-F_i(x))^{1/\varepsilon} $. Hence,
$\varepsilon\log p^\varepsilon(x,M_0) \geq \max_i\log(1-x_if_i(x))$.
Given $c>0$, choose a neighborhood $V_0$ of $M_0$ such that
$\min_ix_if_i(x) \leq 1-e^{-c}$ whenever $x\in V_0$. Then $\varepsilon
\log p^\varepsilon(x,0) \geq -c $ whenever $x\in V_0$.
\end{pf}

\section*{Acknowledgments}
The authors would also like to thank the participants of the
Probability workshop of the University of Neuch\^atel for useful
remarks and comments.

%suskaldyti doi

% imsref loaded by linak, 2013-09-17 16:51:40
%

\printaddresses

\end{document}